# STOCHASTIC EQUIVARIANT COHOMOLOGIES AND CYCLIC COHOMOLOGY

By Rémi Léandre

*Université de Bourgogne*

We give two stochastic diffeologies on the free loop space which allow us to define stochastic equivariant cohomology theories in the Chen–Souriau sense and to establish a link with cyclic cohomology. With the second one, we can establish a stochastic fixed point theorem.

**1. Introduction.** Let us consider a finite-dimensional orientable manifold $M$ of even dimension. Let us suppose that it is endowed with an action of the circle $S^1$, that is, a smooth map of groups $t \to \psi_t$ from the circle into the set of diffeomorphisms of the manifold. We can suppose since the circle is compact that it is an action by isometries. $d/dt\psi_0 = X$ is called the Killing vector field on $M$.

Duistermaat and Heckman [16] and Berline and Vergne [8] have considered an integral of the following type:

$$(1.1) \qquad \int_M \exp[-dX - |X|^2] \wedge \mu$$

and have shown that this integral is equal to an integral over the fixed point set of the circle action, that is, the manifold where the Killing vector field vanishes. We have only, in order to show this localization formula, to suppose that $(d + i_X)\mu = 0$ (this means that $\mu$ is equivariantly closed) and that $\mu$ is of even degree. In order to understand (1.1), let us remark that we can endow $M$ with a Riemannian structure invariant under the circle action. $X$ can be considered alternatively as a vector field or as a 1-form: in (1.1), $dX$ is considered as the exterior derivative of the 1-form $X$ and $|X|^2 = i_X X$ is a scalar. Forms of even degree constitute a commutative algebra. Therefore $\exp[-dX - |X|^2] = \exp[-|X|^2] \sum \frac{(-1)^n}{n!} dX^{\wedge n}$ and in (1.1), we consider the integral of the top degree form.

This leads to the concept of $S^1$-equivariant cohomology of $M$. We consider the set of invariant forms under the circle action on $M$. If a form is invariant under rotation, its Lie derivative for the Killing vector field $(d + i_X)^2 = di_X + i_X d$ is equal to 0. This shows that the equivariant exterior derivative $d + i_X$ defines a complex on the set of forms invariant by rotation. This complex is called the $S^1$-equivariant complex. The main theorem of Jones and Petrack [30] is that the $S^1$-equivariant cohomology is equal to the de Rham cohomology of the fixed point set under the circle action of the manifold.

Let us recall how Jones and Petrack proceed in order to prove this theorem. They remark that the fixed point set of the Killing vector field (the set where $X = 0$) is a manifold. They assume that there is a neighborhood $T$ of the fixed point set which is invariant under rotation and which retracts equivariantly on $\{X = 0\}$. They deduce that the equivariant cohomology of $T$ is equal to the de Rham cohomology of the fixed point set: namely, the equivariant cohomology of the fixed point set is equal to the traditional de Rham cohomology of it, because $X = 0$ on the fixed point set. They remark that the equivariant cohomology of $\{X \neq 0\}$ is trivial: namely, $X$ considered as a 1-form is such that $(d + i_X)X$ is invertible in the algebra of forms invariant under rotation on $\{X \neq 0\}$. They conclude by a Mayer–Vietoris argument: namely, the equivariant cohomology of $T$ is equal to the de Rham cohomology of the fixed point set, and the equivariant cohomology of $T \cap \{X \neq 0\}$ and of $\{X \neq 0\}$ is equal to 0. We get a cover of $M$ by open subset invariant by rotation, such that we can apply the mechanism of the long exact sequence of Mayer and Vietoris.

We are interested in an infinite-dimensional generalization of this work. Namely, in theoretical physics, people consider the free loop space of $M$ of smooth maps $\gamma$ from the circle $S^1$ into $M$. Let us consider a compact spin manifold such that the smooth free loop space $L_\infty(M)$ of smooth maps $\gamma$ from the circle $S_1$ into $M$ is orientable. It carries a natural circle action, and the fixed point set is the manifold itself. The generator of this circle action $X_\infty$, called the canonical Killing vector field, is the vector over a loop which to $s$ associates $\frac{d\gamma}{ds}$. Namely, the fiber in $\gamma$ of the tangent space of $L_\infty(M)$ coincides with the space of smooth sections of the bundle on the circle $\gamma^* T(M)$ where $T(M)$ is the tangent bundle of $M$. It consists of smooth maps from the circle $s \to X(s)$ into $T(M)$ such that $X(s)$ belongs to the fiber of the tangent bundle on $\gamma(s)$. $s \to X_\infty(s) = d/ds\, \gamma(s)$ is such a smooth map because we consider the smooth loop space.

Following Atiyah [6], the index of the Dirac operator $D_+$ over $M$ should satisfy

$$(1.2) \qquad \mathrm{Ind}\, D_+ = C \int_{L_\infty(M)} \exp[-dX_\infty - |X_\infty|^2]$$

such that the index theorem over the manifold should be a localization formula in infinite dimension in the manner of Duistermaat and Heckmann [16] or Berline and Vergne [8]. In (1.2), Atiyah considered the $L^2$ metric on the tangent space $\gamma^* T(M)$ of a loop $\gamma$. In particular, we remark thar $|X_\infty|^2$ is nothing else but the energy of the loop $\gamma \int_{S^1} |d/ds\gamma(s)|^2 \, ds$. Bismut [10, 11] pioneered the relation between the equivariant cohomology of the loop space and the index theory by considering the Dirac operator over $M$ tensorized by an auxiliary bundle. He introduced the Bismut–Chern character over the free loop space, which is equivariantly closed, and which is related to the index theorem of the Dirac operator. The Bismut–Chern character is associated to the bundle $\xi_\infty$ on the free loop space deduced from the bundle $\xi$ on $M$ as follows: the fiber of $\xi_\infty$ on $\gamma$ consists of smooth sections of the bundle on $S^1(\gamma^*\xi)$. The reader interested in further developments about this topic by physicists can see the book of Szabo [60] and the references therein.

So the equivariant cohomology of the free loop spaces gives topological invariants. Jones and Petrack [30] show that the equivariant cohomology of the smooth loop space is equal to the cohomology of the manifold, by localization. Getzler, Jones and Petrack [22] introduce iterated integrals and establish a link between the equivariant cohomology of the loop space and the cyclic cohomology. Inspired by this, Getzler [21] defines algebraically a current over the loop space, which by localization gives the index theorem.

Our motivation is to give an analytical meaning to the current of Getzler [21]. For that, we need a measure over the free loop space. We choose the B–H–K measure, which is invariant under rotation [10, 11, 25]. Let us recall quickly the definition of the B–H–K measure. Let $\Delta$ be the Laplace–Beltrami operator on $M$. Let $p_t(x, y)$ be the heat kernel associated to the heat semigroup on $M$. Let $dP_{1,x}$ be the law of the Brownian bridge starting from $x$ and coming back at time 1 to $x$. The B–H–K measure is given by

$$(1.3) \qquad d\mu = \frac{p_1(x, x) \, dx \otimes dP_{1,x}}{\int_M p_1(x, x) \, dx}.$$

In some sense, it is the unique measure on the free loop space which is invariant under rotation [18] and which is constructed from the Brownian bridge measure.

The first remark is that the equivariant cohomology of the free loop space is related to a series of forms of arbitrary degree. Jones and Léandre [29] have introduced a Hilbert tangent space over a random loop, which was given in a preliminary form by Bismut [9]. This allows us to define an $L^p$ theory of forms over the loop space and to show that the Bismut–Chern character belongs to all the $L^p$. By using the integration by parts over the free loop space of Léandre [33, 38], we can establish in the line of Malliavin calculus a Sobolev cohomology theory over the Brownian bridge, show that the Sobolev cohomology groups of the loop space are equal to the Hochschild cohomology groups [35, 37, 41, 43, 45] if the manifold is simply connected, and show therefore that the Sobolev cohomology groups are equal to the cohomology groups of the smooth loop space.

In functional analysis, there are differential calculi that are different from the classical calculus. Let us recall, for instance, what is a Frölicher space [31] (or a space endowed with a smooth structure). A space $M$ endowed with a vector space $F_M$ of maps from $M$ into $\mathbb{R}$ and a set $C_M$ of maps from $\mathbb{R}$ into $M$ is called a Frölicher space if and only if the two following conditions hold:

(i) $f$ belongs to $F_M$ if and only if for all $c$ in $C_M$, $f \circ c$ is a smooth map from $\mathbb{R}$ into $\mathbb{R}$.

(ii) $c$ belongs to $C_M$ if and only if for all $f$ in $F_M$, $f \circ c$ is a smooth map from $\mathbb{R}$ into $\mathbb{R}$.

Let $(M, F_M, C_M)$ and $(M', F'_M, C'_M)$ be two Frölicher spaces. A map $\phi$ from $M$ into $M'$ is said to be smooth if $f' \circ \phi$ belongs to $F_M$ as soon as $f'$ belongs to $F'_M$, or, equivalently, if $\phi \circ c$ belongs to $C'_M$ as soon as $c$ belongs to $C_M$ [31].

There is another calculus, which deals with forms, and which was introduced by Chen and Souriau, which is analogous to Frölicher calculus. Let us recall what a diffeology on a topological space $M$ is. It is constituted of a collection of maps $(\phi_U, U)$ from any open subset $U$ of any $\mathbb{R}^n$ satisfying the following requirements:

(i) If $j : U_1 \to U_2$ is a smooth map from $U_1$ into $U_2$, and if $(\phi_{U_2}, U_2)$ is a plot, $(\phi_{U_2} \circ j, U_1)$ is still a plot called the composite plot.

(ii) The constant map is a plot.

(iii) If $U_1$ and $U_2$ are two open disjoint subsets of $\mathbb{R}^n$ and if $(\phi_{U_1}, U_1)$ and $(\phi_{U_2}, U_2)$ are two plots, the union map $\phi_{U_1 \cup U_2}$ realizes a plot from $U_1 \cup U_2$ into $M$.

This allows Chen and Souriau to define a form. A form $\sigma$ is given by the data of forms $\phi_U^* \sigma$ on $U$ associated to each plot $(\phi_U, U)$. The system of forms over $U$ $\phi_U^* \sigma$ has moreover to satisfy the following requirement: if $(\phi_{U_2} \circ j, U_1)$ is a composite plot, $(\phi_{U_2} \circ j)^* \sigma$ is equal to $j^* \phi_{U_2}^* \sigma$.

The exterior derivative $d\sigma$ of $\sigma$ is given by the data $d\phi_U^* \sigma$.

The main example of Souriau is the following: let $M$ be a manifold endowed with an equivalence relation $\sim$. We can consider the quotient space $\tilde{M}$. Let $\pi$ be the projection from $M$ onto $\tilde{M}$. A map $\tilde{\phi}$ from an open subset $U$ of a finite-dimensional linear space is a plot with values in $\tilde{M}$ if, *by definition*, there is a smooth lift $\phi$ from $U$ into $M$ such that $\tilde{\phi} = \pi \circ \phi$.

The ideas of Chen and Souriau lead to another stochastic differential calculus, which deals with forms almost surely defined as in Malliavin calculus, and which is more flexible: it is the stochastic Chen–Souriau calculus (see [13, 26, 59] in the deterministic case). There are many diffeologies which lead to different stochastic de Rham cohomology theories, but in general these stochastic cohomology theories are equal to the de Rham cohomology groups of the smooth loop space or the Hölder loop space [40, 44, 46, 49].

The fact that the stochastic Chen–Souriau calculus is more flexible than the Sobolev calculus allows us in this present work to deal with the stochastic equivariant cohomology of the free loop space, although $X_\infty$ is not defined over the Brownian bridge, because the Brownian loop is not differentiable.

In the first part, we define a poor diffeology, which is very simple, and allows us to define a stochastic equivariant cohomology of the free loop space. There are a few stochastic plots, such that we get many forms, smooth in the Chen–Souriau sense: let us recall, for example, that if a diffeology is included in a second one, a form for the second one is still a form for the first one. We require that the operation $\psi_t$ gotten by rotating a loop is smooth for the considered diffeology, that it transforms a stochastic plot into a stochastic plot. This shows that the set of stochastic forms is invariant under rotation. We establish a link between the equivariant cohomology in this sense and the cyclic cohomology. We show that the stochastic equivariant cohomology of $L(M) - M$ is zero, but we cannot prove by using this diffeology a fixed point theorem, because we cannot produce a retract of an equivariant small contractible neighborhood of the constant loop which is compatible with this diffeology. Namely, we have to produce a retract which satisfies the two requirements: it is smooth for the considered diffeology and commutes with $\psi_t$ for all $t$ in $S^1$.

For that reason, we consider in the second part a richer diffeology, but more artificial than the first one, which allows us to produce this retract. This gives a fixed point theorem: the equivariant cohomology with respect to this diffeology is equal to the cohomology of $M$. The relation with cyclic cohomology is performed by using the theory of anticipative Stratonovitch integrals of Léandre [35] over the loop space.

Let us recall the previous work in order to define $X_\infty$: in [38], the stochastic Killing vector field is defined as an antisymmetric operator of order 1, which is densely defined, therefore closable. But this construction does not work for forms. In [36] the interior product by the stochastic Killing vector field is defined as a fermionic Hida distribution [24], but the program failed because the iterated integral does not belong in the domain of this distribution. White noise analysis [24] has defined the derivative of the flat Brownian motion in another way. Léandre [51] has defined the speed of the curved Brownian bridge as a white noise distribution operating on stochastic iterated Chen integrals. Léandre [50] has considered the case of the hypoelliptic bridge; the difference with the work here is that the considerations of Léandre [50] are not intrinsic, because Hörmander's type operator is written under a nonintrinsic form. Moreover, the relation with Léandre [50] and index theory is not clear.

The reader can see the two surveys of Léandre about analysis over loop space and topology [42, 48] and the survey of Albeverio [2] about analysis on loop space and mathematical physics.

**2. Study of the first diffeology.** Let us consider the free loop space of finite energy $L_1(M)$, that is, the set of maps $\gamma$ from the circle $S_1$ into $M$ such that

$$(2.1) \qquad \int_0^1 |d/ds\,\gamma(s)|^2 \, ds < \infty.$$

It is a Hilbert manifold. A deterministic plot $\phi_{\det}$ of dimension $m$ is given by the following data:

  (a) $U$ an open subset of $\mathbb{R}^m$,
  (b) a smooth map $\phi_{\det}$ from $U$ into $L_1(M)$.

The set of all deterministic plots of $L_1(M)$ constitutes a diffeology [13, 26, 59].

DEFINITION 2.1. A deterministic form $\sigma_{\det}$ on $L_1(M)$ is given by the following: to each plot $\phi_{\det}$, we associate a smooth form $\sigma_U = \phi_{\det}^* \sigma_{\det}$ over $U$. Moreover, the set of finite-dimensional forms $\phi_{\det}^* \sigma_{\det}$ satisfies the following property: if $j : U_1 \to U_2$ is a smooth map and if $\phi_{2,\det}$ is a plot with parameter space $U_2$, we can consider the plot $\phi_{1,\det} = \phi_{2,\det} \circ j$. Then

$$(2.2) \qquad \phi_{1,\det}^* \sigma_{\det} = j^* \phi_{2,\det}^* \sigma_{\det}.$$

REMARK. Since a smooth function from $U$ into $\mathbb{R}$ is smooth if and only if its restriction to each smooth path in $U$ is smooth, it is equivalent to saying that a functional is smooth in the Chen–Souriau sense or is smooth in the Frölicher sense. We take as $C_{L_1(M)}$ the space of smooth curves from $\mathbb{R}$ into $L_1(M)$ and as $C_{L_1(M)}$ the space of maps $f$ from $L_1(M)$ into $\mathbb{R}$ such that $f \circ c$ is smooth for all $c$ in $C_{L_1(M)}$. But $L_1(M)$ is a Hilbert manifold. A curve $c$ from $R$ into a Hilbert space $H$ is smooth if for all $h$ in $H$, $\langle h, c \rangle$ is smooth from $R$ into $R$ ([31], Theorem 2.14).

REMARK. Our notion of form is an adaptation in our situation of the notion of functional smooth in the Gateaux sense on the loop space.

A Frechet smooth form gives a deterministic form in this sense. We can consider $n$ vector fields $X_i$ on $U$. $\phi_{\det}^* \sigma_{\det} = \sigma(\phi_{\det})(D_{X_1}\phi_{\det}, \ldots, D_{X_n}\phi_{\det})$ because since $\phi_{\det}$ is Frechet-smooth, $D_{X_i}\phi_{\det}$ realizes naturally an element of the tangent bundle of $L_1(M)$ in $\phi_{\det}$. A tangent vector $s \to X_s$ on a loop belonging to $L_1(M)$ can be written $s \to \tau_s H_s$ where $\tau_s$ is the parallel transport on the loop. Moreover, $s \to H_s$ is of finite energy. Therefore the fiber over a loop of the tangent bundle of the loop space is a Hilbert space. We can define the cotangent bundle as usual and the $n$-exterior power of the tangent bundle. We get a Hilbert bundle $\Lambda^n(L_1(M))$ of $n$-form over $L_1(M)$ and an $n$-form is a smooth section of this bundle.

On the free loop space, there is a natural circle action $\psi_t : \gamma \to \{s \to \gamma(t + s)\}$, which is a smooth transformation of the finite energy loop space. Its generator is called $X_{\infty,\det}$ and is not a vector field over $L_1(M)$: $X_{\infty,\det}(\gamma)(s) = d/ds\,\gamma(s)$. (In a more convenient way, we should look to the smooth loop space in order to speak of a smooth circle action, which is endowed with the strucure of a Frechet manifold, and replace the previous considerations by the Frechet topology on the smooth loop space. In the sequel, we should replace the polygonal approximation by approximation by convolution.) We can consider forms in the present sense weaker than an ordinary form over $L_\infty(M)$, the smooth free loop space endowed with the Frechet topology. As before, a traditional form over $L_\infty(M)$ is a form in this sense. Namely, a vector over a smooth loop $\gamma$ coincides to a smooth section over $\gamma$ of the tangent bundle of $M$, such that an $n$-form coincides with an $n$-antilinear distribution, which depends smoothly on $\gamma$ in $L_\infty(M)$. A tangent vector over a loop $s \to X_s$ is written $s \to \tau_s H_s$ where $s \to H_s$ is smooth such that $\tau_1 H_1 = H_0$. The tangent space of $\gamma$ is a Frechet space. The cotangent space of $\gamma$ coincides with the dual of the Frechet space which gives the tangent space of $\gamma$. We denote it $T_\gamma^*$. We can consider the bundle of $n$ forms $\Lambda^n(L_\infty(M))$ of $L_\infty(M)$. An $n$-form can be seen as the set of alternated continuous forms on $T_\gamma$. Since the tangent bundle of $L_\infty(M)$ is locally trivialized and since $\gamma \to \tau_1$ is smooth for the Frechet topology on $L_\infty(M)$, an $n$-form can be seen as a smooth section of $\Lambda^n(L_\infty(M))$ (see [31] for analogous discussions). The previous considerations are easier to see on the based loop space $L_{x,\infty}(M)$ of loops starting from $x$ and arriving at $x$. A tangent vector $X(\cdot)$ on a loop $\gamma$ can be seen as $X(s) = \tau(s)H(s)$ where $\tau(s)$ is the parallel transport from $\gamma(0)$ to $\gamma(s)$ along the loop $\gamma$ and $H(\cdot)$ is a smooth path in $T_{\gamma(0)}(M)$ such that $H(0) = H(1) = 0$. The tangent bundle of the smooth based loop space is therefore trivial. $T_\gamma^*$ can therefore be realized as a fixed space of distributions, endowed with its dual topology. A 1-form is a smooth application in the Frechet sense in this space of distributions. A natural extension can be done for the definition of Frechet-smooth $n$-form on the smooth based loop space. In the sequel, we will use the notion of forms smooth in the Chen–Souriau sense weaker than the traditional definition of forms smooth in the Frechet sense, because it is more consistent with the framework of this work. We define the exterior derivative of a deterministic form $\sigma_{\det}$ by the set of relations $d(\phi_{\det}^*\sigma_{\det}) = \phi^*(d\sigma_{\det})$ for any plot. This checks clearly the relations of Definition 2.1.

In the following, we will use the notion of extended plot $\phi_{\det}^{\text{ext}}(u,t)$ of a plot: $u \in U$; $t \in S_1$:

$$(2.3) \qquad \phi_{\det}^{\text{ext}}(u,t)(s) = \psi_t \phi_{\det}(u)(s).$$

DEFINITION 2.2.   $i_{X_{\infty,\det}}\sigma_{\det}$ is given for a plot $\phi_{\det}$ by

$$(2.4) \qquad \phi_{\det}^{*} i_{X_{\infty,\det}}\sigma_{\det}(u) = i_{\partial/\partial t}\phi_{\det}^{\text{ext}*}\sigma_{\det}(u,0).$$

DEFINITION 2.3.   A form $\sigma_{\det}$ is called invariant under the circle action if for all $t$ frozen

$$(2.5) \qquad \phi_{\det}^{\text{ext},*}\sigma_{\det}(u,t) = \phi_{\det}^{\text{ext},*}\sigma_{\det}(u,0).$$

Let us recall that $\phi_{\det}^{\text{ext},*}\sigma_{\det}(u,t)$ is a form on $U \times S^1$. $\frac{\partial}{\partial t}$ is a vector field on $U \times S^1$. $[i_{\partial/\partial t}\phi_{\det}^{\text{ext}*}\sigma_{\det}(u,0)$ is the form $\phi_{\det}^{\text{ext}*}\sigma_{\det}(u,0)(\frac{\partial}{\partial t},\cdot)$.] If $\sigma_{\det}$ is issued from a form in the classical sense over the Frechet manifold $L_{\infty}(M)$, we have that $\phi_{\det}^{\text{ext}*}\sigma_{\det}(u,0)(\cdot,\frac{\partial}{\partial t}) = \sigma_{\det}(\cdot,\frac{\partial}{\partial t}\psi_0(\phi_{\det}(u)))$. But the quantity $\frac{\partial}{\partial t}\psi_0(\phi_{\det}(u))$ is nothing else than $X_{\infty,\det}(\phi_{\det}(u))$. So Definition 2.2 is consistent.

Over $U \times S^1$, we have a natural circle action. An $n$-form is invariant under the circle action over the free loop space if for all plots, the form $\phi_{\det}^{\text{ext}*}\sigma_{\det}$ is invariant under the circle action on $U \times S^1$. $X_{\infty,\det}$ corresponds to the vector field $\frac{\partial}{\partial t}$ on $U \times S^1$. On $U \times S^1$, the invariant forms under rotation are written $\sigma(u) + \sigma_1(u) \wedge dt$ where $\sigma(u)$ and $\sigma_1(u)$ are forms on $U$ which do not depend on $t$.

If a deterministic form is invariant under rotation, we have, by seeing plots:

$$(2.6) \qquad d(i_{X_{\infty,\det}}\sigma_{\det}) + i_{X_{\infty,\det}}(d\sigma_{\det}) = 0.$$

Namely, we can consider the Lie derivative of the $\phi_{\det}^{\text{ext}*}\sigma_{\det}(u,t)$, which does not depend on $t$, in the $t$ direction. Then we use the formula expressing the Lie derivative $L_X$ along a vector field $X$ in terms of the exterior derivative and the interior product along the vector field: $L_X = di_X + i_X d$. Since $d^2 = 0$ and since a double interior product $i_{X_{\infty,\det}}i_{X_{\infty,\det}}\sigma_{\det} = 0$, we deduce the following theorem:

Let $\Lambda_{\det}^{\text{ev}}$ be the set of formal series of deterministic forms invariant by rotation over $L_{\infty}(M)$ of even degree and let $\Lambda_{\det}^{\text{odd}}$ be the set of formal series of deterministic forms invariant by rotations of odd degrees. We define $\Lambda_{\det}^{2k} = \Lambda_{\det}^{\text{ev}}$ and $\Lambda_{\det}^{2k+1} = \Lambda_{\det}^{\text{odd}}$.

THEOREM 2.4.   $d + i_{X_{\infty,\det}}$ *realizes a complex from* $\Lambda_{\det}^{k}$ *into* $\Lambda_{\det}^{k+1}$ *for all $k$.*

Let $L(M)$ be the continuous free loop space. Let $\Delta$ be the Laplace–Beltrami operator associated to the Riemannian metric over the compact manifold, which is imbedded isometrically in $\mathbb{R}^d$. The heat semigroup has a heat kernel $p_t(x, y)$. Let $dP_{1,x}$ be the law of the Brownian bridge starting from $x$ and coming back at time 1 at $x$. We put (see [10, 11, 25]):

$$(2.7) \qquad d\mu = \frac{p_1(x, x)\, dx \otimes dP_{1,x}}{\int_M p_1(x, x)\, dx}.$$

It is a probability measure over $L(M)$, which is invariant under the natural circle action on the loop space.

DEFINITION 2.5. A stochastic plot of dimension $m$ $\phi_{\mathrm{st}} = (U, \phi_i, \Omega_i)_{i \in \mathbb{N}}$ is given by the following data:

(a) a fixed open subset $U$ of $\mathbb{R}^m$,

(b) a countable measurable partition $\Omega_i$ of $L(M)$,

(c) a family of smooth applications $(u, s, y) \to F_i(u, s, y)$ from $U \times S_1 \times M$ bounded with bounded derivatives of all orders (if we work initially over the finite energy loop space, we suppose only it has finite energy in $s$) and a family of application $r_i : U \to S^1$ constants on each connected component of $U$,

(d) over $\Omega_i$, $\phi_i(u) = \{s \to F_i(u, s, \gamma(s + r_i(u)))\}$ belongs to $L(M)$.

REMARK. The system of $F_i$ tells us how we deform the random loop $\gamma$: we allow to deform the random loop via cylindrical functional as, for instance, $F_{(s, s, \gamma_s)} = \exp_{\gamma_i(s)}[u(\gamma(s) - \gamma_i(s))]$ where $\exp_{\gamma_i(s)}[v]$ is the Riemannian exponential centered in $\gamma_i(s)$ and where $v$ is a vector in $T_{\gamma_i(s)}$. $(\gamma(s) - \gamma_i(s)) = w$ is the unique vector such that $\exp_{\gamma_i(s)}[w] = \gamma(s)$ for a smooth loop $\gamma_i$ close from $\gamma$.) In the third part, we will give an extended way to deform the loop $\gamma$: we refer to [43] for a way of deforming the loop $\gamma$ in all the class of semimartingales.

We remark that $\psi_t \phi_{\mathrm{st}}$ is still a stochastic plot. It is given by $s \to F_i(u, s + t, \gamma(s + r_i(u)))$ on $\psi_t \Omega_i$. (We consider the cover $\psi_t \Omega_i$ of the free loop space instead of the cover $\Omega_i$.)

Let $\Omega^N$ be the set of loops such that $\sup_{|s-t|<1/N} d(\gamma(s), \gamma(t)) < r$ where $d$ is the Riemannian distance over $M$ and $r$ is a small positive real number. By considering the partition $\Omega_i \cap \Omega^N$, we can suppose that in Definition 2.5, each $\Omega_i$ is imbedded in an $\Omega^N$. If $\gamma \in \Omega^N$, we denote by $\gamma^N$ its polygonal approximation by broken geodesics, if we work over the finite energy loop space. If we work over the smooth loop space, we regularize the loop $\gamma$ as follows: we consider its convolution in $\mathbb{R}^d$, where the regularizing kernel is a support smaller than $1/N$; we get a loop $\tilde{\gamma}^N$ in $\mathbb{R}^d$ which is never far from $M$.

We consider the projection function $\pi$ from a tubular neighborhood of $M$ (which is supposed imbedded in $\mathbb{R}^d$) into $M$ conveniently extended to $R^d$, and we put $\gamma^N = \pi\tilde{\gamma}^N$. For the stochastic integrals which are considered in this part, the two types of approximations lead to the same result, because we consider nonanticipative Stratonovitch integrals, but in the next part, this will lead to some complications. Let $\phi_{\mathrm{st}} = (U, \phi_i, \Omega_i)$ be a stochastic plot. We define the approximated plot in the deterministic loop space $\phi_{\mathrm{st}}^N$ of length $N$ associated to $\phi_{\mathrm{st}}$ by: if $\Omega_i \subseteq \Omega^N$, $\phi_i^N(u) = \{s \to F_i(u, s, \gamma(s + r_i(u)))^N\}$ over $\phi_i^{-1}\Omega^N$.

Let us remark that $\phi_i^{-1}O^N$ is included into $\phi_i^{-1}\Omega^{N+1}$ and their union is equal to $U$ over $\Omega_i$.

$\phi_{\mathrm{st}}^N$ defines a random plot over $L_1(M)$, or in order to be more rigorous, a stochastic plot over $L_\infty(M)$ the smooth free loop space, if we use the convolution approximation. We can give the notion of stochastic form in a way which is a little bit more sophisticated than in [40]. In the remaining part of this work, we will work always on $L_\infty(M)$.

DEFINITION 2.6.    A stochastic form $\sigma_{\mathrm{st}}$ is given by the following data:

 (i) A deterministic form $\sigma_{\det}$ over $L_\infty(M)$ which is called the skeleton of the deterministic form.

 (ii) To each stochastic plot $\phi_{\mathrm{st}} = (U, \phi_i, \Omega_i)_{i\in\mathbb{N}}$, $\phi_{\mathrm{st}}^{N*}\sigma_{\det}$ tends in probability to a random form over $U$ $\phi_{\mathrm{st}}^*\sigma_{\mathrm{st}}$ for the smooth topology with compact support. This means that over each $\Omega_i$, over each $\phi_i^{-1}\Omega^N$, $\phi_i^{N'*}\sigma_{\det}$ for $N' > N$ tends for this topology to the restriction of $\phi_{\mathrm{st}}^*\sigma_{\mathrm{st}}$ to this open subset of $U$.

REMARK.    $\sigma_{\det}$ is only defined by $\sigma_{\mathrm{st}}$, because we can consider the plot $u \to \{s \to F(u, s)\}$ with values in $L_\infty(M)$. Namely, a smooth function $\phi_U : U \to L_\infty(M)$ can be seen as a smooth function form $U \times S^1$ with values in $M(u, s) \to \phi_U(s)$.

REMARK.    Let $j : U_1 \to U_2$ be a deterministic map from $U_1$ into $U_2$. Let $\phi_{\mathrm{st}}^2 = (U_2, \phi_i^2, \Omega_i)_{i\in N}$ be a stochastic plot and let $\phi_{\mathrm{st}}^1 = (U_1, \phi_i^2 \circ j, \Omega_i)_{i\in\mathbb{N}}$ be the composite plot. We get almost surely as random form

$$j^* \phi_{\mathrm{st}}^{2*}\sigma_{\mathrm{st}} = \phi_{\mathrm{st}}^{1*}\sigma_{\mathrm{st}}. \tag{2.8}$$

Moreover, let $\phi_{\mathrm{st}}^1 = (U, \phi_i^1, \Omega_i^1)_{i\in\mathbb{N}}$ and $\phi_{\mathrm{st}}^j = (U, \phi_j^2, \Omega_j^2)_{j\in\mathbb{N}}$ be two stochastic plots. Let us suppose that there exists a random transformation $\Psi$ from some $\Omega_i^1$ into some $\Omega_j^2$ such that $\phi_j^2 = \phi_i^2 \circ \Psi$. Then, almost surely,

$$\phi_j^{1*}\sigma_{\mathrm{st}} = \phi_i^{2*}\sigma_{\mathrm{st}} \circ \Psi \tag{2.9}$$

as random form over $U$.

The two properties (2.8) and (2.9) constitute the basic properties of a random form in our previous works [40, 46].

A random form $\sigma_{\mathrm{st}}$ is said to be invariant by rotation if its skeleton $\sigma_{\mathrm{det}}$ is invariant by rotation. By the second property of Definition 2.6, if $\sigma_{\mathrm{st}}$ is a random form with skeleton $\sigma_{\mathrm{det}}$, $d\sigma_{\mathrm{st}}$ is still a random form and its skeleton is $d\sigma_{\mathrm{det}}$. Namely, $d\sigma_{\mathrm{st}}$ satisfies clearly the requirement of Definition 2.6 because $\phi_{\mathrm{st}}^{N*}d\sigma_{\mathrm{det}} = d\phi_{\mathrm{st}}^{N*}\sigma_{\mathrm{det}}$. Therefore the stochastic exterior derivative defines a complex over the set of random form. It is not the same for $i_{X_\infty}\sigma_{\mathrm{st}}$, where $X_\infty$ is the stochastic Killing vector field which formally generates the circle action over the free loop space.

DEFINITION 2.7. A stochastic form $\sigma_{\mathrm{st}}$ with skeleton $\sigma_{\mathrm{det}}$ has an interior product $i_{X_\infty}$ if $i_{X_{\infty,\mathrm{det}}}\sigma_{\mathrm{det}}$ defines a stochastic form. This stochastic form is called $i_{X_\infty}\sigma_{\mathrm{st}}$.

REMARK. Let us see what these definitions mean with some examples. If $F_{\mathrm{det}} = \int_{S^1}|d/ds\gamma(s)|^2\,ds$, $F_{\mathrm{det}}$ does not define a stochastic functional. If $\omega$ is a one form on $M$, $F_{\mathrm{det}} = \int_{S^1}\langle\omega(\gamma(s)), d\gamma(s)\rangle$ defines a stochastic functional. Let $e_s$ be the evaluation map $\gamma(\cdot) \to \gamma(s)$. $e_s^*\omega = \sigma_{\mathrm{det}}$ defines a stochastic form, but this stochastic form does not admit an interior product by the stochastic Killing vector field, unlike the stochastic form $\int_{S^1}e_s^*\omega\,ds$.

Clearly, if $i_{X_\infty}\sigma_{\mathrm{st}}$ exists, $i_{X_\infty}i_{X_\infty}\sigma_{\mathrm{st}} = 0$. If $i_{X_\infty}\sigma_{\mathrm{st}}$ exists, $di_{X_\infty}\sigma_{\mathrm{st}}$ exists, but it is not clear that $i_{X_\infty}d\sigma_{\mathrm{st}}$ exists.

DEFINITION 2.8. $\Lambda_{\mathrm{st}}^{2k}$ is the set of formal series of even stochastic forms which are invariant under rotation and which admit an interior product by $X_\infty$. $\Lambda_{\mathrm{st}}^{2k+1}$ is the set of formal series of odd stochastic forms which are invariant under rotation and which admit an interior product by $X_\infty$.

If a form of given degree invariant under rotation admits an interior product by $X_\infty$, its skeleton $\sigma_{\mathrm{det}}$ satisfies

$$(2.10) \qquad d(i_{X_{\infty,\mathrm{det}}}\sigma_{\mathrm{det}}) + i_{X_{\infty,\mathrm{det}}}(d\sigma_{\mathrm{det}}) = 0.$$

Therefore $i_{X_\infty}d\sigma_{\mathrm{st}}$ exists and we have the relation

$$(2.11) \qquad d(i_{X_\infty}\sigma_{\mathrm{st}}) + i_{X_\infty}(d\sigma_{\mathrm{st}}) = 0.$$

We therefore get:

THEOREM 2.9. $d + i_{X_\infty}$ defines a complex from $\Lambda_{\mathrm{st}}^k$ to $\Lambda_{\mathrm{st}}^{k+1}$. It is called the stochastic equivariant complex.

Its cohomology groups $H_{st}^{odd}$ and $H_{st}^{ev}$ are called the stochastic equivariant cohomology groups with respect to this diffeology.

Let us give an example. Let $\tilde{\omega}^n$ be an element of the algebraic tensor product $\Omega(M) \otimes \Omega.(M)^{\otimes n-1}$ where $\Omega(M)$ denotes the space of smooth forms over $M$ and $\Omega.(M)$ denotes the space of smooth forms of degree not equal to 0. If $\tilde{\omega}^n = \omega_1 \otimes \omega_2 \otimes \cdots \otimes \omega_n$, we define its degree by

$$(2.12) \qquad \deg \tilde{\omega}^n = \deg \omega_1 + \sum_{i=1}^{n}(\deg(\omega_i) - 1).$$

In the sequel, we will denote by $C^{2k}$ the sum of finite sums of elementary tensor products of even degree. We get an analogous definition of $C^{2k+1}$.

Let $\tilde{\omega}^n = \omega_1 \otimes \omega_2 \otimes \cdots \otimes \omega_n$. Let us introduce the Hochschild boundary:

$$
\begin{aligned}
(2.13) \qquad b(\tilde{\omega}^n) = {} & \sum_{i=1}^{n-1} \varepsilon_i \omega_1 \otimes \cdots \otimes \omega_i \wedge \omega_{i+1} \cdots \otimes \omega_n \\
& + \varepsilon_n \omega_n \wedge \omega_1 \otimes \omega_2 \cdots \otimes \omega_{n-1} \\
& + \sum_{i=1}^{n} \varepsilon_i' \omega_1 \otimes \cdots \otimes \omega_{i-1} \otimes d\omega_i \otimes \cdots \otimes \omega_n.
\end{aligned}
$$

The signs $\varepsilon_i$ and $\varepsilon_i'$ are given in [22]. The Hochschild boundary increases the degree by one unit.

Let us define the cyclic boundary operator $B$ of Connes [14, 22, 53]:

$$(2.14) \qquad B(\tilde{\omega}^n) = \sum_{i=1}^{n} \varepsilon_i'' 1 \otimes \omega_i \cdots \otimes \omega_n \otimes \omega_1 \otimes \omega_2 \otimes \cdots \otimes \omega_{i-1}.$$

The signs $\varepsilon_i''$ are given in [22]. $B$ decreases the degree by one unit.

Let us recall (see [22]) that $b + B$ realizes a complex from $C^k$ into $C^{k+1}$ called the cyclic complex. Let us introduce the equivariant stochastic integral $\Sigma \tilde{\omega}^n$ defined by

$$(2.15) \quad \Sigma \tilde{\omega}^n = \int_0^1 e_s^* \omega_1 \wedge \int_{s < s_2 < s_3 < \cdots < s_n < s+1} e_{s_2}^* \omega_2(d\gamma_{s_2}, \cdot) \wedge \cdots \wedge e_{s_n}^* \omega_n(d\gamma_{s_n}, \cdot)$$

where $e_s$ denotes the evaluation map $\gamma \to \gamma(s)$. Let us clarify what we mean by $e_s^* \omega(d\gamma_s, \cdot)$. It is nothing else than $i_{X_{\infty,det}} e_s^* \omega(\cdot) \, ds$, if we work over $L_\infty(M)$. $\Sigma \tilde{\omega}^n$ realizes a smooth form over $L_\infty(M)$. The Hochschild boundary corresponds to the exterior product of the iterated integral. The Connes operator corresponds to the interior product of the iterated integral by the canonical deterministic Killing vector field $X_{\infty,det}$. In general, Chen forms are obtained when we remove the first integral $\int_0^1$ in (2.15) and we take $s = 0$. These nonequivariant Chen forms will lead to stochastic forms which do not admit an interior product by the stochastic Killing vector field, unlike the equivariant one in (2.15).

We get from the rules of approximation of nonanticipative Stratonovitch integrals the following theorem:

THEOREM 2.10. $\Sigma\tilde{\omega}^n$ defines a random form $\Sigma_{st}\tilde{\omega}^n$.

Let us explain Theorem 2.10. Let us consider $\tilde{\omega} = \omega_1 \otimes \omega_2 \otimes \cdots \otimes \omega_n$, where the degree of $\omega_i$ is equal to $r_i$. We extend the form $\omega_i$ by forms with bounded derivatives of all orders over $R^d$. Let $\phi_i(u) = \{s \to F_i(u, s, \gamma(s + r_i(u)))\}$ be a stochastic plot and let $\phi_i^N$ be its regularization by convolution. Let $X_{i_j}$ be some vector fields over $U$ where $i_j \in [1, r_i - 1]$ if $i > 1$ and $1_j \in [1, r_1]$. Modulo some antisymmetrization which is due to the fact that we take some exterior product in the iterated integrals we get

$$\phi_i^{N*}\Sigma\tilde{\omega}(X_{1_1}, \ldots, X_{1_{r_1}}, X_{2_1}, \ldots, X_{2_{r_2-1}}, \ldots, X_{n_1}, \ldots, X_{n_{r_n-1}})$$

$$= \int_0^1 \langle \omega_1((\phi_i^N)(s_1)), D_{X_{1_1}}\phi_i^N(s_1), \ldots, D_{X_{1_{r_1}}}\phi_i^N(s_1) \rangle \, ds_1$$

(2.16)
$$\times \int_{s_1 < s_2 < \cdots < s_n < s_1 + 1} \prod \langle \omega_j(\phi_i^N(s_j)), d/ds_j \phi_i^N(s_j),$$

$$D_{X_{j_1}}\phi_i^N(s_j), \ldots, D_{X_{j_{r_j-1}}}\phi_i^N(s_j) \rangle \, ds_j$$

$$+ c.c.$$

The counterterm appears when we shuffle the $X_{j,k}$, because we consider a wedge product of forms in the definition of the iterated integrals. The stochastic integral converges in probability for the smooth topology on $U$ to the limit nonanticipative Stratonovitch integral, where we replace formally in (2.16) $\phi_i^N(s_j)$ by $\phi_i(s_j)$ and $d/ds_j\phi_i^N(s_j)$ by the Stratonovitch differential $d_{s_j}\phi_i(s_j)$ (see [27]). After this formal manipulation, we get a random smooth form on $U$.

Let us justify this limiting argument: let us consider a semimartingale $s \to Y_s$ with values in $M$ and some 1-form $\omega_i$ on $M$. We define iteratively $I_{n,t} = \int_0^t I_{n-1,s} \langle \omega_n(Y_s), dY_s \rangle$ with $I_{1,t} = \int_0^t \langle \omega_1(Y_s), dY_s \rangle$. We can replace $Y_t$ by $Y_t^N$. We get random variables $I_{n,t}^N$. By induction, we can show that $I_{n,t}^N$ tends in probability to $I_{n,t}$ (see [27]).

Let us explain the formulas by a simpler example. Let us suppose that $\tilde{\omega} = \omega_1 \otimes \omega_2$ where $\omega_1$ is a 1-form and $\omega_2$ is a 2-form. We extend them to bounded forms with bounded derivatives of all orders over $R^d$. Let $\phi_i(u) = \{s \to F_i(u, s, \gamma(s + r_i(u)))\}$ and let $\phi_i^N$ be the approximated plot. Let $X$ and $Y$ be two canonical vector fields over $\mathbb{R}^m$. We have

$$\phi_i^{N*}\Sigma\tilde{\omega}(X, Y) = \int_0^1 \langle \omega_1(\phi_i^N(s)), \partial_X \phi_i^N(s) \rangle \, ds$$

$$(2.17) \quad \times \int_s^{s+1} \langle \omega_2(\phi_i^N(s)), d/ds\, \phi_i^N(s), \partial_Y \phi_i^N(s) \rangle \, ds$$

$$- \int_0^1 \langle \omega_1(\phi_i^N(s)), \partial_Y \phi_i^N(s) \rangle \, ds$$

$$\times \int_s^{s+1} \langle \omega_2(\phi_i^N(s)), d/ds\, \phi_i^N(s), \partial_X \phi_i^N(s) \rangle \, ds.$$

The counterterm in (2.16) is here very simple to write. These stochastic integrals converge in probability for the smooth topology in $U$ over the limit nonanticipative Stratonovitch integral. The ordinary integral in $ds$ gives an equivariant form. We produce by this procedure a map between the stochastic equivariant cohomology and the cyclic cohomology, which extends this classical correspondence on the smooth loop space (see [22]).

Let $M$ be the fixed point set of $L(M)$ under the circle action. Let $T_\varepsilon(M) = \{\gamma : \sup_{s,t} d(\gamma(s), \gamma(t)) < \varepsilon\}$ for a small real nuimber $\varepsilon$. It is an equivariant neighborhood of $M$. Let $h(s,t) = d^2(\gamma(s), \gamma(t))$ if $\gamma(s)$ and $\gamma(t)$ are close and equal to 1 if $\gamma(s)$ and $\gamma(t)$ are far. We introduce a function $g(h(s,t)) = 1$, if $h(s,t) < r_1$, which behaves as $(r_2 - h(s,t))^{-k}$ for a big $k$ if $h(s,t) \to r_{2-}$ and is equal to infinity if $h(s,t) > r_2$. We suppose $g$ is smooth over $[0, r_2[$ and larger than 1. Let $f$ be a smooth function from $[1, \infty[$ into $[0, 1]$ equal to 1 at 1 and with compact support. We consider the functional

$$(2.18) \quad H^{r_1, r_2}(\gamma) = 1 - f\left( \int_0^1 \int_0^1 g(h(s,t)) \, ds \, dt \right).$$

We define $O_\varepsilon = L(M) - \overline{T}_\varepsilon(M)$. It is an equivariant open subset of $L(M)$. We find $\varepsilon$ and $\varepsilon'$ with $\varepsilon' > \varepsilon$ such that $H^{r_1, r_2}$ with support in $O_\varepsilon$ is equal to 1 in $O_{\varepsilon'}$. *In the sequel, in order to get this property, we will work on the Hölder loop space.* It is not a problem, because the Brownian bridge is almost surely Hölder. If the first property is satisfied, we say $H^{r_1, r_2}$ satisfies property H. The cutoff functional $H^{r_1, r_2}$ is Frechet smooth for the Hölder topology, and a fortiori smooth for the smooth topology (see [46]).

We say that a stochastic form belongs to $\Lambda_{\mathrm{st}}(O_\varepsilon)$ if for all mollifier functions $H^{r_1, r_2}$ satisfying the property H the form $H^{r_1, r_2}\sigma_{\mathrm{det}}$ is a stochastic form over the full loop space $L(M)$. This defines a skeleton $\sigma_{\mathrm{det}}$ over $L_\infty(M) \cap O_\varepsilon$ such that for all mollifier functions $H^{r_1, r_2}\sigma_{\mathrm{st}}$ determines a stochastic form over the free loop space.

Since $O_\varepsilon$ is invariant by rotation, we can repeat the considerations of before for stochastic forms over $O_\varepsilon$, define an interior product by the stochastic Killing vector field and define an equivariant exterior derivative. We get two stochastic cohomology groups $H_{\mathrm{st}}^{\mathrm{odd}}(O_\varepsilon)$ and $H_{\mathrm{st}}^{\mathrm{ev}}(O_\varepsilon)$. The end of this part is devoted to showing the following proposition:

PROPOSITION 2.11.   $H_{\mathrm{st}}^{\mathrm{ev}}(O_\varepsilon) = H_{\mathrm{st}}^{\mathrm{odd}}(O_\varepsilon) = 0.$

The main difficulty to repeat the argument of Jones and Petrack [30] is that the Brownian loop is of infinite energy.

Let us remark that over $O_\varepsilon$, $\int_0^1 |d/ds\gamma(s)|^2 ds > a > 0$ for a given constant $a$. Let us consider the cover of $]0, \infty[$ by the intervals $]\frac{a}{8(k+1)}, \frac{a}{4k}[$ where $\frac{a}{0} = \infty$. Let us consider a partition of unity $f^k$ associated to this cover. Let us introduce a function $g$ from $\mathbb{R}$ into $[0, 1]$ equal to 1 over $\mathbb{R}^-$ and equal to zero over $[a/2, \infty[$. We imbed $M$ into $\mathbb{R}^d$, and we denote by $\langle \cdot, \cdot \rangle_x$ the Riemannian tensor in $x$ in $\mathbb{R}^d$ which extends the Riemannian metric over $M$. We put

$$(2.19) \quad F^{N,k}(\gamma) = f^k\left(\int_0^1 dt \int_0^1 \langle d/ds\psi_t(\gamma^N)(s), d/ds\psi_t(\gamma)(s) \rangle_{\psi_t(\gamma)(s)} ds\right).$$

$F^{N,k}$ is invariant under rotation. We put

$$(2.20) \quad G^N(\gamma) = \prod_{i=1}^{N-1} g\left(\int_0^1 dt \int_0^1 \langle d/ds\psi_t(\gamma^i)(s), d/ds\psi_t(\gamma)(s) \rangle_{\psi_t(\gamma)(s)} ds\right).$$

$G^N$ is invariant under rotation. $G^N$ and $F^{N,k}$ define clearly stochastic functionals with respect to the diffeology, by replacing $d/ds\psi_t\gamma(s)\, ds$ by the stochastic integral $d_s\psi_t\phi_i(s)$. Moreover, over $L_\infty(M)$, we have

$$(2.21) \qquad\qquad \infty > \sum_{N,k} G^N F^{N,k} = \Xi > 0.$$

The sum is in fact finite. Let $O^{N,k}$ be the open subset:

$$(2.22) \quad \begin{aligned} &\left\{\gamma : \int_0^1 \int_0^1 ds\, dt \langle d/ds\psi_t(\gamma^N)(s),\right. \\ &\qquad\qquad \left. d/ds\psi_t(\gamma)(s) \rangle_{\psi_t(\gamma)(s)} \in \left] \frac{a}{8(k+1)}, \frac{a}{4k} \right[ \right\}. \end{aligned}$$

$\tilde{F}^{N,k} = \frac{G^N F^{N,k}}{\Xi}$ constitutes a partition of unity associated to the cover $O^{N,k}$ of $L_\infty(M) \cap O_\varepsilon$. Moreover $\tilde{F}^{N,k}$ constitutes a smooth stochastic functional with respect to the previous diffeology. It constitutes a stochastic partition of unity of $O_\varepsilon$. That is, if we consider a stochastic plot $\phi_{\mathrm{st}} = (U, \phi_i, \Omega_i)_{i \in \mathbb{N}}$, we have almost surely as smooth function over $U$ for all integers $M$:

$$(2.23) \qquad\qquad \sum_{N,k} \tilde{F}^{N,k}(\phi_i^M(u)) = 1$$

over each $\Omega_i$ and the sum is almost surely finite.

In the sequel, we put $(N, k) = \alpha$. Over the set of indices, there is a natural order. We consider a set of indices $(\alpha_1 < \cdots < \alpha_n) = I_n$. $O^{I_n} = O^{\alpha_1, \ldots, \alpha_n} = \bigcap O^{\alpha_i}$.

We say that a stochastic form is defined over $O^{I_n}$ if for all smooth function $h^{\bar{j}}$ with compact support in $]\frac{a}{8(k_j+1)}, \frac{2}{4k_j}[$, the form $\prod h^{\bar{j}}(\int_0^1 \int_0^1 ds\, dt <d/ds\psi_t(\gamma^{N_j})(s), d_s\psi_t(\gamma)(s) >_{\psi_t(\gamma)(s)})\sigma_{\mathrm{st}}$ defines a stochastic form over $L(M)$.

Since $O^{I_n}$ is invariant by rotation, we can define a stochastic equivariant complex over $O^{I_n}$. We get equivariant stochastic cohomology groups called $H^{\mathrm{ev}}_{\mathrm{st}}(I_n)$ and $H^{\mathrm{odd}}_{\mathrm{st}}(I_n)$.

LEMMA 2.12.   *If $n \neq 0$, $H_{\mathrm{st}}(I_n) = 0$.*

PROOF.   Let us consider $N_1$ such that $\alpha_1 = (N_1, k_1)$. We consider the stochastic equivariant 1-form:

$$(2.24) \qquad \alpha^{N_1} = \int_0^1 ds \int_0^1 dt \langle d/ds\psi_t(\gamma^{N_1})(s), \cdot\rangle_{\psi_t(\gamma)(s)}.$$

For a stochastic plot $\phi = (U, \phi_i, \Omega_i)_{i\in N}$, we get if $X$ is a vector field over $U$:

$$(2.25) \qquad \begin{aligned} &\phi^*\alpha^{N_1}(X) \\ &\qquad = \int_0^1 dt \int_0^1 \langle d/ds\psi_t(\phi(u))^{N_1}(s), \partial_X\psi_t(\phi(u))(s)\rangle_{\psi_t(\phi(u))(s)}\, ds. \end{aligned}$$

We have

$$(2.26) \qquad i_{X_\infty}\alpha^{N_1} = \int_0^1 dt \int_0^1 \langle d/ds\psi_t(\gamma^{N_1})(s), d_s\psi_t(\gamma)(s)\rangle_{\psi_t(\gamma)(s)}.$$

This interior product is therefore larger than $\frac{a}{8(k_1+1)}$. Therefore,

$$(2.27) \qquad (d + i_{X_\infty})\alpha^{N_1} = (i_{X_\infty}\alpha^{N_1})\left(1 + \frac{d\alpha^{N_1}}{i_{X_\infty}\alpha^{N_1}}\right).$$

$i_{X_\infty}\alpha^{N_1}$ is a functional in our weaker sense strictly larger than $\frac{a}{8(k_1+1)}$ over $O^{I_n}$. We can define $\frac{1}{i_{X_\infty}\alpha^{N_1}}$. It is therefore still a functional in our sense in $O^{I_n}$. We see that the stochastic forms invariant under rotation over $O^{I_n}$ having an interior product by the stochastic Killing vector field constitute an algebra for the stochastic wedge product. Moreover,

$$(2.28)\ (d+i_{X_\infty})(\sigma_{\mathrm{st}} \wedge \sigma'_{\mathrm{st}}) = (d+i_{X_\infty})\sigma_{\mathrm{st}} \wedge \sigma'_{\mathrm{st}} + (-1)^{\deg\sigma_{\mathrm{st}}}\sigma_{\mathrm{st}} \wedge (d+i_{X_\infty})\sigma'_{\mathrm{st}}$$

by pulling back this formula through a plot and approximating the plot. Therefore, $(d+i_{X_\infty})\alpha^{N_1}$ has an inverse over the algebra of stochastic forms over $O^{I_n}$: it is given by the formula

$$(2.29) \qquad \beta^{N_1} = (i_{X_\infty}\alpha^{N_1})^{-1}\left(\sum(-1)^j \frac{(d\alpha^{N_1})^j}{(i_{X_\infty}\alpha^{N_1})^j}\right).$$

Moreover,

$$(2.30) \qquad (d + i_{X_\infty})\beta^{N_1} = 0$$

by (2.27) and (2.28), because $\alpha^{N_1}$ is invariant under rotation and because

$$(2.31) \qquad \beta^{N_1} \wedge (d + i_{X_\infty})\alpha^{N_1} = 1.$$

It remains to remark that

$$(2.32) \qquad (d + i_{X_\infty})(\sigma \wedge \sigma') = ((d + i_{X_\infty})\sigma) \wedge \sigma' + \text{sign}\,\sigma \wedge ((d + i_{X_\infty})\sigma').$$

Let us choose a stochastic form $\sigma_{\text{st}}$ over $O^{I_n}$ such that $(d + i_{X_\infty})\sigma_{\text{st}} = 0$. We deduce that

$$(2.33) \qquad \sigma_{\text{st}} = (d + i_{X_\infty})(\alpha^{N_1} \wedge \beta^{N_1} \wedge \sigma_{\text{st}}).$$

Therefore the result holds (see [30] for a proof of the same result in the deterministic case). $\quad\square$

Let $I_n - \alpha_j = I_n^j$. A stochastic form associated to $I_n^j$ clearly defines a stochastic form associated to $I_n$. There is a bigraduation over the form associated to $I_n$: the degree in the sense of Theorem 2.4 and the length of $I_n$. We call the space of associated stochastic form $\Lambda^r_{\text{st},I_n}$. We deduce a bicomplex:

(a) the equivariant stochastic derivative $d + i_{X_\infty}$ which transforms $\Lambda^r_{\text{st},I_n}$ into $\Lambda^{r+1}_{\text{st},I_n}$,

(b) the Cech complex:

$$(2.34) \qquad (\delta\sigma_{\text{st}})_{I_n} = \sum (-1)^j \sigma_{\text{st},I_n - \alpha_j}$$

where we take in (2.34) the restriction of $\sigma_{\text{st},I_n-\alpha_j}$ to $O^{I_n}$.

These two complexes commute. We can conclude now:

PROOF OF PROPOSITION 2.11. Let us recall, when there is a first quadrant bicomplex $(d, \delta)$ with complexes commuting $(d\delta = \delta d)$ $\Lambda^{p,q}$, we can consider its total complex $d + \text{sign}\,\delta = d_{\text{tot}}$ which operates on $E^n = \sum_{p+q=n} \Lambda^{p,q}$. $d_{\text{tot}}$ goes from $E^n$ into $E^{n+1}$ and is such that $d_{\text{tot}}^2 = 0$. There are some total cohomology groups which are related to the total complex $\text{Ker}_{d_{\text{tot}}} E^{n+1}/\text{Im}_{d_{\text{tot}}} E^n$. There is an algebraic procedure in order to approximate the total cohomology groups called spectral sequence. We approximate iteratively the total cohomology by a sequence of cohomology groups associated to some complexes $d_r$ which applies $\Lambda^{p,q}$ to $\Lambda^{p+r,q-r+1}$. More precisely, at the further step, we start from the cohomology groups of the previous step. At a first step, we consider $d$ alone, and we get $E_1^{p,q} = H_d(\Lambda^{p,q})$, the cohomology groups deduced from $d$. At a second step, we consider the differential $\delta$ alone, and we apply that to the spaces obtained from the second step, and we get $E_2^{p,q} = H_\delta(H_d(\Lambda^{p,q}))$. This approximations procedure converges when we iterate to the total cohomology of the total complex (see [12]). We would like to apply the spectral sequence formalism in our context, mimicking the classical proof of Bott and Tu ([12], pages 166–167) that the Cech cohomology of a manifold is equal to the de Rham cohomology of the manifold. For that, we replace the bicomplex $(\delta, d)$ by the bicomplex $(\delta, d + i_{X_\infty})$. Following Warner ([61], page 202), if $\delta\sigma_{st} = 0$, then $\sigma_{st} = \delta\sigma_{st,1}$, because of the existence of partition of unity *invariant under rotation* on the open subset $O_{\bar{\varepsilon}}$ of the Hölder loop space given previously. We have

$$(2.35) \qquad (\sigma_{st,1})_{I_n} = \sum_\alpha \tilde{F}^\alpha(\sigma_{st})_{\alpha,\alpha_1,\dots,\alpha_n}.$$

Namely, $\tilde{F}^\alpha$ defines a stochastic functional invariant under rotation with support included in $O^\alpha$ such that $\tilde{F}^\alpha(\sigma_{st})_{\alpha,\alpha_1,\dots,\alpha_n}$ is a stochastic form on $O^{I_n}$ invariant under rotation. Since $\tilde{F}^\alpha$ is invariant under rotation, we get a form which is invariant under rotation. We can repeat the proof of Bott and Tu ([12], pages 166–167). In the spectral sequence associated to the bicomplex $(\delta, d + i_{X_\infty})$, the first terms $E_1^{p,q}$ are equal to 0, except in the first column where we find $\Lambda_{st}^p(O_\varepsilon)$. The second terms $E_2^{p,q}$ are all 0 except for the first column where we find the cohomology groups for $d + i_{X_\infty}$ of $\Lambda_{st}^p(O_\varepsilon)$. Therefore, the spectral sequence degenerates after the second order, and we find that the stochastic equivariant cohomology groups of $O_\varepsilon$ are equal to the total cohomology groups of the bicomplex, because the higher derivatives $d^r$ are trivially 0 and at each step of the spectral sequence, the cohomology groups considered remain invariant.

We invert the role of the two complexes. Since the stochastic cohomology for the equivariant stochastic exterior derivative over $O^{I_n}$ is equal to zero, all the terms $E_1^{p,q}$ are equal to 0. This shows that the total stochastic cohomology groups of the total bicomplex are equal to 0, because the first cohomology groups at each step of the spectral sequence are 0, and because the cohomology groups at step $r + 1$ of the spectral sequence are deduced from the cohomology groups at step $r$ of the spectral sequence, supposed inductively equal to 0.

Therefore the result holds, because we have computed the total cohomology groups of the bicomplex in two ways, which are equal.  $\square$

We can state an analogous proposition:

PROPOSITION 2.13. *Let $\varepsilon' > \varepsilon > 0$. Then the stochastic equivariant cohomology groups of $O_\varepsilon \cap T_{\varepsilon'}$ are equal to 0.*

**3. Study of the second diffeology.** Let $T_\varepsilon$ be the equivariant open subset defined by $\{\gamma : \sup_{s,t} d(\gamma(s), \gamma(t)) < \varepsilon\}$. We would like to show that its equivariant stochastic cohomology groups are equal to the de Rham cohomology groups of $M$, if $\varepsilon$ is small enough. We need a retraction map from $T_\varepsilon$ to $M$ which commutes with the natural circle action. Let $r \in [0, 1]$. We choose

$$(3.1) \qquad F(r, \gamma(s), \gamma(t)) = \exp_{\gamma(s)}[r(\gamma(t) - \gamma(s))]$$

conveniently extended over $\mathbb{R}^d$ by a functional with values in $\mathbb{R}^d$ with bounded derivatives of all orders. $\gamma(t) - \gamma(s)$ is the vector over $\gamma(s)$ of the unique geodesic joining $\gamma(s)$ to $\gamma(t)$ if $\gamma(t)$ is closed from $\gamma(s)$. This gives a retraction map $F(r)$ from the loop $t \to \gamma(t)$ to the constant loop $t \to \gamma(s)$. Moreover, $F(r)$ transforms a plot into a plot; that is, $F(r)$ is smooth for the considered diffeology. But we do not have $\psi_t F(r) = F(r)\psi_t$, because we contract the small loop $t \to \gamma(t)$ into the constant loop $t \to \gamma(s)$. We request a retract $H(r)$ which is smooth for the considered diffeology and which commutes with the circle action. $F(r)$ is not equivariant under the natural circle action because we choose the time $s$. We average under the natural circle action: we get $\int_{S_1} F(r, \gamma(s), \gamma(t)) \, ds$ which is not far of $M$ if $\gamma \in T(\varepsilon)$ and we look at the projection map $\pi$ conveniently extended over $R^d$ on $M$: $\pi(\int_{S_1} F(r, \gamma(s), \gamma(t)) \, ds) = H(r, \gamma)(t)$. The map $H(r, \cdot)$ commutes with the natural circle action. But $H(r)$ does not transform a plot into a plot. This leads to the introduction of a new stochastic diffeology.

DEFINITION 3.1. A stochastic plot $\phi_{\text{st}} = (U, \phi_i, \Omega_i)_{i \in \mathbb{N}}$ is given by the following data:

(i) any finite sequence of deterministic integers $j$,

(ii) a deterministic open subset $U$ of $\mathbb{R}^m$,

(iii) a countable measurable partition $\Omega_i$ of $L(M)$,

(iv) two applications $F_i^j$ from $U \times (\mathbb{R}^d)^{n_i^j} \times \mathbb{R}^d$ into $\mathbb{R}^d$ and $h_i^j$ from $U \times \mathbb{R}^d \times \mathbb{R}^d$ smooth with bounded derivatives of all orders and an application $r_i^j$ from $U$ into $S^1$ which is constant on the connected component of $U$,

(v) let us denote $H_{h_i^j, F_i^j}(u)(s)$ the quantity $\phi_i(u)(s) = h_i^j(\int_{(S_1)^n} F_i^j(u, \gamma(s_1),$ $\ldots, \gamma(s_n), \gamma(s + r_i^j(u))) \, ds_1 \cdots ds_n)$. The iteration $H_{h_i^1, F_i^1} \circ \cdots \circ H_{h_i^n, F_i^n}(u)(\cdot)$ belongs to $L(M)$ over $\Omega_i$.

The main remark is the following: if $\phi_{\text{st}} = (U, \phi_i, \Omega_i)$ is a plot, $(r, u) \to H(r, \phi_{\text{st}}(u))$ is still a plot indexed by $U \times [0, 1]$. This stochastic diffeology is compatible with the restriction map. If $\phi_{\text{st}} = (U, \phi_i, \Omega_i)_{i \in N}$ is a plot with respect to this diffeology, we get an extended plot $\phi_{\text{st}}^{\text{ext}}$ from $U \times [0, 1]$ into $L(M)$ by putting

$$(3.2) \qquad (u, r) \to \{s \to H(r, \phi_{\text{st}}(u))(s)\}$$

which contracts the stochastic plot $\phi_{\mathrm{st}}$ with values in $T(\varepsilon)$ into a plot with values in $M$. This says that the retraction map is smooth for this new diffeology.

We can repeat the consideration of Section 2 to study the stochastic equivariant cohomology associated to this diffeology.

In particular, we get:

PROPOSITION 3.2.    $H_{\mathrm{st}}^{\mathrm{ev}}(O_\varepsilon) = H_{\mathrm{st}}^{\mathrm{odd}}(O_\varepsilon) = 0$.

We also get:

PROPOSITION 3.3.    If $\varepsilon' > \varepsilon > 0$, the stochastic equivariant cohomology of $O_\varepsilon \cap T_{\varepsilon'}$ is equal to 0.

In order to show this theorem, we do as in the previous section. There is a small difficulty which appears, because in $\int_0^1 \langle d/ds\psi_t(\gamma^N)(s), d/ds\psi_t(\gamma(s))\rangle_{\psi_t(\gamma(s))}\, ds$, there are some anticipative Stratonovitch stochastic integrals which appear. We replace this expression by $\int_0^1 \langle d/ds\psi_t(\gamma^N)(s), d/ds\psi_t(\gamma(s))\rangle_{\psi_t(\gamma^N)(s)}\, ds$ and we integrate by parts in order to remove the stochastic integral, if we replace $d/ds\psi_t(\gamma(s))\, ds$ by the anticipative Stratonovitch differential $d_s\psi_t(\phi_{\mathrm{st}})$ where $\phi_{\mathrm{st}}$ is a plot.

PROPOSITION 3.4.    The stochastic equivariant cohomology groups of $T_\varepsilon$ are equal to the cohomology groups of $M$ if $\varepsilon > 0$ is small enough.

PROOF.    Let $H(r, \cdot)$ be the application $\gamma \to \{s \to H(r, \gamma)(s)\}$. It commutes over the smooth loop space to the circle action. Therefore over the smooth loop space, we have

(3.3)                    $X_{\infty,\mathrm{det}}(H(r,\gamma)) = DH(r,\gamma)X_{\infty,\mathrm{det}}$.

Namely, $\psi_t \circ H(r) = H(r) \circ \psi_t$ and we differentiate this formula at $t = 0$.

We denote by $X_r$ the vector $\frac{\partial}{\partial r}H(r, \cdot)$. This realizes a Frechet vector field on $L_\infty(M)$; more precisely, $X_r(\gamma)$ belongs to $T_{H^r(\gamma)}$. We use the retraction Cartan formula for a deterministic form $\sigma_{\mathrm{det}}$ over $T(\varepsilon)$. We get

$$\sigma_{\mathrm{det}} = H^{0*}\sigma_{\mathrm{det}} + d\int_0^1 H^{r,*}i_{X_r}\sigma_{\mathrm{det}}\, dr + \int H^{r*}i_{X_r}(d + i_{X_{\infty,\mathrm{det}}})\sigma_{\mathrm{det}}\, dr$$

(3.4)

$$- \int_0^1 H^{r,*}i_{X_r}i_{X_{\infty,\mathrm{det}}}\sigma_{\mathrm{det}}\, dr.$$

If $\sigma_{\mathrm{det}}$ were a traditional form over $L_\infty(M)$, this formula would be nothing else than the integrated formula which expresses the Lie derivatives along a flow in terms of the exterior derivative and the interior product along the vector field of the flow of the considered deterministic form $d/dr \, H^{r*}\sigma_{\det} = dH^{r*}i_{X_r}\sigma_{\det} + H^{r*}i_{X_r}\,d\sigma_{\det}$. But, here, we have to take care, because we consider a weaker notion of form. So, we have to look at this formula through a plot, and consider the extended retracted plot. We apply this formula to the form associated to the finite-dimensional form which is given by the extended plot. By (3.3),

$$(3.5) \qquad H^{r,*}i_{X_{\infty,\det}}\sigma_{\det} = i_{X_{\infty,\det}}H^{r,*}\sigma_{\det}.$$

Therefore, if $\sigma_{\det}$ is equivariantly closed, then we have the equivariant retraction formula:

$$(3.6) \qquad \sigma_{\det} = H^{0,*}\sigma_{\det} + (d + i_{X_{\infty,\det}})\int_0^1 H^{r,*}i_{X_r}\sigma_{\det}\,dr.$$

This formula is still true for the stochastic form. Namely, if $\phi_{\mathrm{st}}$ is the stochastic plot, we have an augmented stochastic plot $(r, u) \to H(r, \phi_{\mathrm{st}}(u))$ called $\phi_{\mathrm{st}}^{\mathrm{aug}}$. We can define $\phi_{\mathrm{st}}^{\mathrm{aug},*}\sigma_{\mathrm{st}}$, and its approximation $(r, u) \to (\phi_{r,t}^{\mathrm{aug}}(r, u))^N$ (see [50] for similar considerations). $H(r)^*(\sigma_{\mathrm{st}})$ is defined by taking the plot $u \to \phi_{\mathrm{st}}^{\mathrm{aug}}(r, u)$ where $r$ is frozen. It admits an interior product by the stochastic Killing vector field, and the approximating formula (3.6) when we pull-back $\sigma_{\mathrm{st}}$ by the approximating plots goes to the limit. Therefore, if $(d + i_{X_\infty})\sigma_{\mathrm{st}} = 0$, we have

$$(3.7) \qquad \sigma_{\mathrm{st}} = H(0)^*\sigma_{\mathrm{st}} + (d + i_{X_\infty})\int_0^1 H(r)^*i_{X_r}\sigma_{\mathrm{st}}\,dr.$$

$H(0)^*\sigma_{\det}$ is a stochastic form over $M$, therefore a deterministic smooth form. It remains to show that $H(r)^*i_{X_r}\sigma_{\mathrm{st}}$ is a stochastic form having an interior product by the Killing vector field. But we can use the fact that $i_{X_\infty}H(r)^*i_{X_r}\sigma_{\mathrm{st}} = -H(r)^*i_{X_r}i_{X_\infty}\sigma_{\mathrm{st}}$ and the fact that $\sigma_{\mathrm{st}}$ admits an interior product by the stochastic Killing vector field in order to show this statement. This proves the proposition. $\qquad\square$

By using Propositions 3.2, 3.3 and 3.4, we can show a stochastic fixed point theorem.

**THEOREM 3.5.** *The stochastic equivariant cohomology groups with respect to the stochastic diffeology of Definition 3.3 are equal to the de Rham cohomology groups of $M$.*

**PROOF.** Let $0 < \varepsilon < \varepsilon'$. We have a cover of $L(M)$ by $T(\varepsilon)$ and $O(\varepsilon)$. These subsets are invariant under rotation. We have a partition of unity associated to this cover [see (2.17)] for the Hölder topology, which is invariant under rotation, which is therefore smooth for the Frechet topology over the smooth loop space, and which provides therefore a partition of unity associated to our cover, invariant under rotation, for the stochastic Chen–Souriau calculus, because a functional Frechet-smooth on the Hölder free loop space realizes clearly a functional smooth in the stochastic Chen–Souriau sense. We can produce a Mayer–Vietoris long exact sequence for the stochastic equivariant cohomology (see [12], pages 22–23). This Mayer–Vietoris argument says that we have an exact sequence:

$$(3.8) \qquad 0 \to \Lambda_{st}(L(M)) \to \Lambda_{st}(T_{\varepsilon'}) \oplus \Lambda_{st}(O_{\varepsilon}) \to \Lambda_{st}(T_{\varepsilon'} \cap O_{\varepsilon}) \to 0$$

for stochastic forms invariant under rotation because the mollifier $H^{r_1,r_2}$ is invariant by rotation. The first map gives the restriction of $\sigma_{st}$ to $T_{\varepsilon'}$ and $O_{\varepsilon}$. There is a partition of unity Frechet smooth for the Hölder topology associated to the cover $T_{\varepsilon'}$ and $O_{\varepsilon}$ of $L_{\infty}(M)$. The functions $\rho_{T_{\varepsilon'}}$ and $\rho_{O_{\varepsilon}}$ associated to this partition of unity are invariant under rotation. $(-\rho_{O_{\varepsilon}}\sigma_{st}, \rho_{T_{\varepsilon'}}\sigma_{st})$ realize an element of $\Lambda_{st}(T_{\varepsilon'}) \oplus \Lambda_{st}(O_{\varepsilon})$ which projects on $\sigma_{st}$ which belongs to $\Lambda_{st}(T_{\varepsilon} \cap O_{\varepsilon})$. From this short exact sequence, we deduce a long exact sequence in cohomology (see [12]). We use the complex $d + i_{X_{\infty}}$, which is compatible with the maps of (3.8). This long sequence in cohomology arises by an abstract argument. Propositions 3.2, 3.3 and 3.4 show the result. Namely, $H_{st}(T_{\varepsilon'} \cap O_{\varepsilon}) = 0 = H_{st}(O_{\varepsilon})$ and $H_{st}(T_{\varepsilon'}) = H(M)$. We deduce that $H_{st}(L(M)) = H(M)$ by the Mayer–Vietoris long sequence in cohomology:

$$(3.9) \qquad \begin{aligned} \cdots &\to H_{st}^{\cdot}(L(M)) \to H_{st}^{\cdot}(T_{\varepsilon'}) \oplus H_{st}^{\cdot}(O_{\varepsilon}) \\ &\to H_{st}^{\cdot}(T_{\varepsilon'} \cap O_{\varepsilon}) \to H_{st}^{\cdot+1}(L(M)) \to \cdots, \end{aligned}$$

where the image of a map in (3.9) is equal to the kernel of the map which follows. $\qquad \square$

In order to show that an iterated integral in the manner of (2.15) defines a stochastic form with respect to this diffeology, we have to study the approximation of anticipative Stratonovitch integrals by convolution. It is a refinement of the theory of Léandre [35].

We work over the based path space, that is, the space of continuous paths starting from $x$ endowed with the Brownian motion measure. Let us recall what is the Sobolev Nualart–Pardoux calculus of Léandre [35, 37, 45]. The tangent space of a path $\gamma$ is the set of map $s \to \tau_s H_s$ where $\tau_s$ is the parallel transport along the path $\gamma$ and where $s \to H_s$ is a path in $T_x(M)$ of finite energy. We take as Hilbert norm

$$(3.10) \qquad \|X.\|_{\gamma}^2 = \int_0^1 |d/ds H_s|^2 \, ds.$$

If $H.$ is deterministic, we have an integration by parts formula [9, 15, 33]:

$$E[\langle dF, X \rangle] = E[F \operatorname{div} X], \tag{3.11}$$

where $F$ is a cylindrical functional. This allows us to define the notion of $H$-derivative:

$$\langle dF, X \rangle = \int_0^1 \langle k(s), d/ds H_s \rangle \, ds. \tag{3.12}$$

We can iterate this notion of stochastic derivative, by using the connection $\nabla$ on the path space:

$$\nabla \tau . H. = \tau . DH.. \tag{3.13}$$

If $d_\nabla^r F$ is defined, we put

$$\begin{aligned} (3.14) \qquad & d_\nabla^{r+1} F(X_1, \ldots, X_{r+1}) \\ & = \langle d(d_\nabla^r(X_1, \ldots, X_r)), X_{r+1} \rangle - \sum d_\nabla^r F(X_1, \ldots, \nabla_{X_{r+1}} X_i, \ldots, X_r). \end{aligned}$$

$d_\nabla^r F$ is an $r$-cotensor and is defined by a kernel:

$$\begin{aligned} (3.15) \qquad & d_\nabla^r F(X_1, \ldots, X_r) \\ & = \int \langle k(s_1, \ldots, s_r), d/ds H_{s_1}, \ldots, d/ds H_{s_r} \rangle \, ds_1 \cdots ds_r \end{aligned}$$

and we put as curved Sobolev norm:

$$\|F\|_r^p = E\left[ \left( \int |k(s_1, \ldots, s_r)|^2 \, ds_1 \cdots ds_r \right)^{p/2} \right]^{1/p} \tag{3.16}$$

(see [33, 38]).

We define the Nualart–Pardoux Sobolev norms for $s \to H(s)$ with values in $\mathbb{R}^d$. We consider the kernels of $d_\nabla^r H(s)$ called $H(s, s_1, \ldots, s_r)$. We suppose that outside the diagonals

$$\begin{aligned} (3.17) \qquad & \|H(s, s_1, \ldots, s_2) - H(s', s_1', \ldots, s_r')\|_{L^p} \\ & \leq C_{p,r}(H)\left( \sum \sqrt{|s_i - s_i'|} + \sqrt{|s - s'|} \right) \end{aligned}$$

and we suppose that for all $s, s_1, \ldots, s_r$

$$\|H(s, s_1, \ldots, s_r)\|_{L^p} \leq C_{p,r}'(H) < \infty. \tag{3.18}$$

The smallest quantities $C_{p,r}(H)$ and $C_{p,r}'(H)$ constitute the system of Nualart–Pardoux norms of the process $H(\cdot)$.

We imbed the manifold into $R^d$. We suppose that $h$ and $F$ are bounded functionals from $\mathbb{R}^d$ into $\mathbb{R}^d$ [$F$ from $(\mathbb{R}^d)^{n+1}$ into $\mathbb{R}^d$] with bounded derivatives of all orders. Let us introduce $Y(s)$ which is a finite iteration of operations of the type $h(\int_{[0,1]^n} F(\gamma(s_1), \ldots, \gamma(s_n), \gamma(s + s_0)) \, ds_1 \cdots ds_n)$ and its approximation by convolution:

$$(3.19) \qquad Y^N(s) = \int g_N(s - u) Y(u) \, du.$$

We will choose the regularizing function $g_N$ later.

We remark that $Y(s)$ satisfies the Nualart–Pardoux conditions by the following lemma:

LEMMA 3.6.   *Let us suppose that $s \to Y(s)$ satisfies the Nualart–Pardoux conditions. Then the random process $s \to h(\int_{[0,1]^n} F(Y(s_1), \ldots, Y(s_n), Y(s + s_0)) \, ds_1 \cdots ds_n)$ satisfies the Nualart–Pardoux conditions.*

PROOF.   The proof comes from the proof of Lemme A.2 of [35] and from the fact $(s_1, \ldots, s_n, s) \to (Y(s_1), \ldots, Y(s_n), Y(s))$ satisfies the Nualart–Pardoux conditions $(s_1, \ldots, s_n, s) \in [0, 1]^{n+1}$ included. (We have a natural extension to this case of the notion of Nualart–Pardoux conditions.)   □

By an integration by parts, and using a primitive $g_N^1$ of $g_N$, we get that

$$(3.20) \qquad Y^N(s) = \int g_N^1(s - u) \, d_u Y(u)$$

such that

$$(3.21) \qquad d/ds \, Y^N(s) = \int g_N(s - u) \, d_u Y(u).$$

We have the following lemma:

We choose a regularizing function $g_N$ such that $\frac{g_N}{2N}$ is equal to 1 over $[-1/N + 1/N^k, 1/N - 1/N^k]$ for a big $k$ and which takes its values in $[0, 1]$, and which is equal to 0 outside $[-1/N; 1/N]$, such that the Nualart–Pardoux norms of $\int_{-1/N}^{-1/N+1/N^k} g_N(s - t) \, d_t Y_t$ are bounded by $N^{-j}$ for a big $j$ as well as the Nualart–Pardoux norms of $\int_{1/N-1/N^k}^{1/N} g_N(s - t) \, d_t Y_t$. The sum of these two terms is called $\delta_N(s)$. The integral

$$(3.22) \qquad \int \langle H^N(s), \delta^N(s) \rangle \, ds$$

goes to 0 by the same considerations as below. So in order to simplify the notation, we can replace $dY^N$ by

$$(3.23) \qquad 2N \int_{-1/N}^{1/N} d_u Y(u - s) = 2N \int_{s-1/N}^{s+1/N} d_u Y(u).$$

Therefore we have to study the behavior of

$$(3.24) \qquad 2N \int_0^t \int_{s-1/N}^{s+1/N} \langle I^N(s) \, ds, d_v Y(v) \rangle$$
$$= 2N \int_{-1/N}^{1/N} dv \int_0^t \langle H^N(s+v), dY(s) \rangle.$$

In order to compute $d_u \psi(u)$ where we consider only one iteration we use the chain rule:

$$(3.25) \qquad d_u \psi(u) = h'\left( \int_{[0,1]^n} F(\gamma(s_1), \ldots, \gamma(s_n), \gamma(s)) \, ds_1 \cdots ds_n \right)$$
$$\times \int_{[0,1]^n} \langle dF(\gamma(s_1), \ldots, \gamma(s_n), \gamma(u)), d\gamma(u) \rangle$$
$$= \langle A(u), d\gamma(u) \rangle.$$

The same result is true when we consider a finite number of iterations as in Definition 3.1, where $A(u)$ checks all the Nualart–Pardoux conditions, $u$ included. So we recognize in (3.23) the anticipative Stratonovitch integal:

$$(3.26) \qquad 2N \int_{-1/N}^{1/N} dv \int_0^t \langle H^N(s+v), A(s) \, d\gamma(s) \rangle.$$

We can write $d\gamma(s) = \tau_s dB(s)$ where $dB(s)$ is a flat Brownian motion. The anticipative Stratonovitch integral $\int_0^t \langle H^N(s+v), A(s) \, d\gamma(s) \rangle$ is equal to the anticipative Stratonovitch integral for the flat Brownian motion $B(s)$

$$(3.27) \qquad \int_0^t \langle H^N(s+v), A(s)\tau(s) \, dB(s) \rangle = \int_0^t \langle \tau(s)^{-1} A_s^t H^N(s+v), dB(s) \rangle.$$

The system of Nualart–Pardoux Sobolev norms for the curved Brownian motion is equivalent to the system of Nualart–Pardoux norms for the flat Brownian motion (see [35, 37, 45]). Let $u_s$ be a process with values in $R^d$. Let $u_s(v)$ be the kernel of its flat derivatives. This means that if $h$ is an element of the Cameron–Martin space of the flat Brownian motion, $D_h v_s = \int_0^1 \langle v_s(v), d/dv h(v) \rangle \, dv$. See $\delta_t(u)$ the Skorohod integral in time $t$ (see [55]). Then the Stratonovitch integral

$$(3.28) \qquad \int_0^t \langle u_s, dB(s) \rangle = \delta_t(u) + \frac{1}{2} \int_0^t \left( \lim_{v \to s-} u_s(v) + \lim_{v \to s+} u_s(v) \right) ds$$

(see [55], page 567). (We do as if we were in $R$ in order to simplify the notation.)

Moreover, we can estimate the $L^p$ norm of $\delta_t(u)$ in terms of the $L^p$ norms of $u_s(v)$ (see [54], page 158) (i.e., in terms of the flat first-order Sobolev norms of $s \to u_s$).

Lemma 3.7. $\int_0^t \langle \tau_s^{-1} A_s^t \frac{1}{2}(H^N(s+v) + H^N(s-v)), dB_s \rangle$ tends when $v \to 0$ to the anticipative Stratonovitch integral $\int_0^t \langle \tau_s^{-1} A_s^t H_s^N(s), dB_s \rangle$ in all the $L^p$.

Proof.    Let $\alpha_s(v)$ be the process $s \to \tau_s^{-1} A_s^t H^N(s+v)$. $\alpha_s(v)$ tends to $\alpha_s(0)$ in all the first-order Sobolev space for the flat Brownian motion. Therefore $\delta_s(\alpha.(v + \alpha.(-v)))$ tends to $\delta_s(0)$ in all the $L^p$. Moreover, if $v > 0$,

$$
\begin{aligned}
(3.29) \quad E\bigg[\bigg|\int_0^t \frac{1}{4}\bigg( &\lim_{t \to s-}(\alpha_s(v)(t) + \alpha_s(-v)(t)) \\
&+ \lim_{t \to s+}(\alpha_s(v)(t) + \alpha_s(-v)(t))\bigg) \\
&- \frac{1}{2}\bigg(\lim_{t \to s+} \alpha_s(0)(t) + \lim_{t \to s-} \alpha_s(0)(t)\bigg)\bigg|^p ds\bigg] = 0.
\end{aligned}
$$

Namely, when $v > 0$, $\lim_{t \to s+} \alpha_s(v)(t) = \lim_{t \to s+} \alpha_{s+v}(t)$, $\lim_{t \to s-} \alpha_s(v)(t) = \lim_{t \to s-} \alpha_{s+v}(t)$. Moreover, $\lim_{t \to s+} \alpha_s(-v)(t) = \lim_{t \to s+} \alpha_{s-v}(t)$ and $\lim_{t \to s-} \alpha_s(-v)(t) = \lim_{t \to s-} \alpha_{s-v}(t)$. Since the Nualart–Pardoux conditions are checked outside the diagonals, we see that $\lim_{v \to 0+} \lim_{t \to s+} \alpha_{s+v}(t) = \lim_{v \to 0+} \lim_{t \to s-} \alpha_{s+v}(t) = \lim_{t \to s-} \alpha_s(t)$ and $\lim_{v \to 0+} \lim_{t \to s+} \alpha_{s-v}(t) = \lim_{v \to 0+} \lim_{t \to s-} \alpha_{s-v}(t) = \lim_{t \to s+} \alpha_s(t)$. □

From Lemma 3.7, we deduce that $2N \int_{-1/N}^{1/N} \langle \int_0^t \tau_s^{-1} A_s^t \frac{1}{2}(H(s+v) + H(s-v)), dB(s)\rangle$ tends in all the $L^p$ to $\int_0^t \langle H(s), dB(s)\rangle$.

We would like to get the same theorem for iterated integrals. We choose

$$
(3.30) \qquad \psi_k^j(s) = h_k^j\bigg(\int_{[0,1]^k} F_k^j(\gamma(s_1), \ldots, \gamma(s_n), \gamma(s))\, ds_1 \cdots ds_n\bigg),
$$

where $h_k^j$ and $k^j$ are smooth with bounded derivatives of all orders. By using the composite of the $\psi_k^j$, in $j$ as in Definition 3.1, we deduce an element $Y_k(s)$. We choose a function $F_k$ from $R^d$ into $R^d$ which is smooth with bounded derivatives of all orders. We define inductively

$$
(3.31) \qquad I^{k+1,N}(t) = \int_0^t I^{k,N}(s) \langle F_k(Y_{k+1}^N)(s), dY_{k+1}^N(s)\rangle.
$$

In order to study the convergence of $I^{N+1}(t)$, it is enough to study the convergence of $\frac{1}{2}(\xi(v) + \xi(-v))$ where $\xi(v) = \int_0^t \frac{1}{2}(I^{k,N}(s+v)) \langle \tau_s^{-1} A_s^{k+1,t} F_{k+1}(Y_{k+1}^N)(s+v), dB(s)\rangle$. The flat derivative of $\int_0^t I_{s+v}^{k,N} \langle \tau_s^{-1} A_s^{k+1,t} F_k(Y_{k+1}^N)(s+v), dB(s)\rangle$ is defined by taking formally the derivative under the sign $\int$ (see [35]).

Let $h(s)$ be a process. Let $h(s, s_1, \ldots, s_n)$ be the kernels of its flat derivatives. Let $K$ be a subset of $(1, \ldots, n)$ and let $\varepsilon_K$ be a collection of sign $\varepsilon_j$ associated to the element of $K$. $\varepsilon_K^c$ is the collection of opposite sign. We denote by $h_{K,\varepsilon_K}(s, u_j)_{j \in K^c}$ the limit of $h(s, u_1, \ldots, u_n)$ when $u_j \to s_{\varepsilon_j}$ for $j$ in $K$. If the Nualart–Pardoux conditions are checked, these expressions exist.

We introduce the following hypothesis of recurrence.

HYPOTHESIS $H^1(k)$. *For an iterated integral of length smaller than $k$, if $v > 0$, uniformly in $N$ and in $s \in [0,1]^m$ when $v \to 0$,*

$$
\begin{aligned}
(3.32) \quad E\Bigg[\Bigg|\sum_{K,\varepsilon_K} \int \prod_{j \in K^c} du_j \tfrac{1}{2} &(I^{k,N}_{K,\varepsilon_K}(s+v)(u_j)_{j \in K^c} \\
&+ I^{k,N}_{K,\varepsilon_K}(s-v)(u_j)_{j \in K^c} \\
&- I^{k,N}_{K,\varepsilon_K}(s)(u_j)_{j \in K^c}\Bigg|^p\Bigg] \to 0.
\end{aligned}
$$

In order to show this property, we have to take the derivatives of $I^{k+1,N}(t)$. The only problem is when we do not take the derivative of $dB_s$. This leads to some anticipative Stratonovitch integrals, which can be treated by (3.28). If $H^1(k)$ is checked, $H^1(k+1)$ is checked. $H^1(0)$ is clearly checked.

This shows us that in the approximation procedure, we can replace the integral between $-1/N$ and $1/N$ by a constant. We do the hypothesis of recurrence where $I^k(t)$ denotes the Stratonovitch integral and where we do no approximation:

HYPOTHESIS $H^2(k)$. *For all $n$, we have when $N \to \infty$ uniformly in $s$*

$$
(3.33) \quad E\Bigg[\Bigg|\int \prod_{j \in K^c} du_j \sum_{K,\varepsilon_K} I^{k,N}_{K,\varepsilon_K}(s)(u_j)_{j \in K^c} - I^k_{K,\varepsilon_K}(s)(u_j)_{j \in K^c}\Bigg|^p\Bigg] \to 0.
$$

If $H^2(k)$ is checked, $H^2(k+1)$ is checked. Namely, the only problem is when we do not take the derivative of the last $dB_t$ in $I^{k+1,N}(t)$ which leads to the study of some anticipative Stratonovitch integrals, which can be treated by (3.28). It is checked for $k = 0$.

We deduce:

PROPOSITION 3.8. *$I^N(t)$ converges in all the flat Sobolev spaces to the anticipative Stratonovitch integral*

$$
(3.34) \quad I(t) = \int_{0 < s_1 < s_2 < \cdots < s_n < t} \langle F_1(Y(s_1)), dY_1(s_1)\rangle \cdots \langle F_n(Y_n(s_n)), dY_n(s_n)\rangle.
$$

Let us suppose that $\psi_k(s)$ depends smoothly on a finite-dimensional parameter $u$:

$$
(3.35) \quad Y^j_k(u, s) = h^j_k\Bigg(u, \int_{[0,1]^n} F^k_j(u, \gamma(s_1), \ldots, \gamma(s_n), \gamma(s)) \, ds_1 \cdots ds_n\Bigg)
$$

and let us consider the composite $\psi_k(u, s)$ as in Definition 3.1. The corresponding approximated integral has a smooth version in $u$, and each of its derivatives in a given $u$ converges in all the flat Sobolev spaces to the corresponding iterated anticipative Stratonovitch integral where we take formally the derivatives in $u$ of the terms which appear in this stochastic integral. But we can estimate the $L^p$ norm of these derivatives over the curved based loop space in terms of the flat Sobolev norms of these derivatives over the flat model, by using the tools of quasi-sure analysis (see [1]). We deduce that the derivative in the parameter space of the approximated integrals converges in all the $L^p$ over the curved Brownian bridge to the derivative in the parameter space of the Stratonovitch integral over the Brownian bridge. By the Sobolev imbedding theorem, we deduce that the approximating integral over the based loop space converges for the smooth topology with compact support in the parameter space in $L^2$.

Let $\tilde\omega^n = \omega_1 \otimes \omega_2 \otimes \cdots \otimes \omega_n$ be an element of $\Omega(M) \otimes \Omega.(M)^{n-1}$. We extend the differential forms $\omega_i$ in smooth forms over $R^d$ bounded with bounded derivatives of all orders. Proposition 3.8 allows us to state the following theorem:

THEOREM 3.9. $\Sigma \tilde\omega^n$ defines a stochastic form with respect to the diffeology of Definition 3.1.

## REFERENCES

[1] AIRAULT, H. and MALLIAVIN, P. (1991). *Quasi-Sure Analysis.* Univ. Paris VI.

[2] ALBEVERIO, S. (1996). Loop groups, random gauge fields, Chern–Simons models, strings: Some recent mathematical developments. In *Espaces de Lacets* (R. Léandre, S. Paycha and T. Wuerzbacher, eds.) 5–34. Univ. Strasbourg.

[3] ALBEVERIO, S., DALETSKII, A. and KONDRATIEV, Y. (2000). Stochastic analysis on product manifolds: Dirichlet operators on differential forms. *J. Funct. Anal.* **176** 280–316. MR1784417

[4] ALBEVERIO, S., DALETSKII, A. and LYTVYNOV, Z. (2001). Laplace operator on differential forms over configuration spaces. *J. Geom. Phys.* **37** 14–46. MR1806439

[5] ARAI, A. and MITOMA, I. (1991). De Rham–Hodge–Kodaira decomposition in infinite dimension. *Math. Ann.* **291** 51–73. MR1125007

[6] ATIYAH, M. F. (1985). Circular symmetry and stationary phase approximation. Colloque en l'honneur de Laurent Schwartz (Paris 1984). *Astérisque* **131** 311–323. MR816738

[7] BENDIKOV, A. and LÉANDRE, R. (1999). Regularized Euler–Poincaré number of the infinite dimensional torus. *Infin. Dimens. Anal. Quantum Probab. Relat. Top.* **2** 617–625. MR1810815

[8] BERLINE, N. and VERGNE, M. (1983). Zéros d'un champ de vecteurs et classes caractéristiques équivariantes. *Duke Math. J.* **50** 539–548. MR705039

[9] BISMUT, J. M. (1984). *Large Deviations and the Malliavin Calculus.* Birkhäuser, Basel. MR755001

[10] BISMUT, J. M. (1985). Index theorem and equivariant cohomology on the loop space. *Comm. Math. Phys.* **98** 213–237. MR786574

[11] Bismut, J. M. (1986). Localisation formulas, superconnection and the index theorem for families. *Comm. Math. Phys.* **103** 127–166. MR826861

[12] Bott, R. and Tu, L. W. (1986). *Differential Forms in Algebraic Topology.* Springer, New York. MR658304

[13] Chen, K. T. (1973). Iterated paths integrals of differential forms and loop space homology. *Ann. Math.* **97** 213–237. MR380859

[14] Connes, A. (1988). Entire cyclic cohomology of Banach algebras and characters of $\theta$-summable Fredholm modules. *K-Theory* **1** 519–548. MR953915

[15] Driver, B. (1992). A Cameron–Martin quasi-invariance formula for Brownian motion on compact manifods. *J. Funct. Anal.* **110** 272–376. MR1194990

[16] Duistermaat, J. J. and Heckman, G. J. (1982). On the variation in the cohomology of the symplectic of the reduced phase-space. *Invent. Math.* **69** 259–269. MR674406

[17] Elworthy, K. D. and Li, X. M. (2000). Special Itô maps and an $L^2$ Hodge theory for one forms on path spaces. Preprint. MR1803385

[18] Emery, M. and Léandre, R. (1990). Sur une formule de Bismut. *Séminaire de Probabilités XXIV. Lecture Notes in Math.* **1426** 448–452. Springer, New York. MR1071560

[19] Fang, S. and Franchi, J. (1997). De Rham–Hodge–Kodaira operator on loop groups. *J. Funct. Anal.* **148** 391–407. MR1469347

[20] Frölicher, A. (1982). Smooth structures. *Category Theory. Lecture Notes in Math.* **963** 69–81. Springer, New York. MR682945

[21] Getzler, E. (1988). Cyclic cohomology and the path integral of the Dirac operator. Preprint. MR936774

[22] Getzler, E., Jones, J. D. S. and Petrack, S. (1991). Differential forms on a loop space and the cyclic bar complex. *Topology* **30** 339–371. MR1113683

[23] Goodwillie, T. G. (1985). Cyclic homology, derivations and the free loop space. *Topology* **24** 187–217. MR793184

[24] Hida, T., Kuo, H. H., Potthoff, J. and Streit, L. (1993). *White Noise: An Infinite Dimensional Calculus.* Kluwer, Dordrecht. MR1244577

[25] Hoegh-Krohn, R. (1974). Relativistic quantum statistical mechanics in two dimensional space time. *Comm. Math. Phys.* **38** 195–224. MR366313

[26] Iglésias, P. (1985). Thesis. Univ. Provence.

[27] Ikeda, N. and Watanabe, S. (1981). *Stochastic Differential Equations and Diffusion Processes.* North-Holland, Amsterdam. MR637061

[28] Jones, J. D. S. (1980). Cyclic cohomology and equivariant cohomology. *Invent. Math.* **87** 403–423. MR870737

[29] Jones, J. D. S. and Léandre, R. (1991). $L^p$ Chen forms on loop spaces. In *Stochastic Analysis* (M. Barlow and N. Bingham, eds.) 104–162. Cambridge Univ. Press. MR1166409

[30] Jones, J. D. S. and Petrack, S. (1990). The fixed point theorem in equivariant cohomology. *Trans. Amer. Math. Soc.* **322** 35–49. MR1010411

[31] Kriegl, A. and Michor, P. W. (1997). *The Convenient Setting of Global Analysis.* Amer. Math. Soc., Providence, RI. MR1471480

[32] Kusuoka, S. (1990). De Rham cohomology of Wiener–Riemannian manifold. Preprint. MR1159292

[33] Léandre, R. (1993). Integration by parts and rotationally invariant Sobolev calculus on free loop space. *J. Geom. Phys.* **11** 517–528. MR1230447

[34] LÉANDRE, R. (1995). Stochastic Moore loop space. *Chaos*: *The Interplay Between Stochastic and Deterministic Behaviour. Lecture Notes in Phys.* **457** 479–502. Springer, Berlin. MR1452632

[35] LÉANDRE, R. (1996). Cohomologie de Bismut–Nualart–Pardoux et cohomologie de Hochschild entiere. *Séminaire de Probabilités XXX. Lecture Notes in Math.* **1626** 68–100. Springer, New York. MR1459477

[36] LÉANDRE, R. (1996). The circle as a fermionic distribution. In *Stochastic Analysis and Related Topics V* (H. Korezlioglu, B. Oksendal and A. S. Ustunel, eds.) 233–236. Birkhäuser, Basel. MR1396333

[37] LÉANDRE, R. (1997). Brownian cohomology of an homogeneous manifold. In *New Trends in Stochastic Analysis* (K. D. Elworthy, S. Kusuoka and I. Shigekawa, eds.) 305–347. World Scientific, Singapore. MR1654376

[38] LÉANDRE, R. (1997). Invariant Sobolev calculus on the free loop space. *Acta Appl. Math.* **46** 267–350. MR1440476

[39] LÉANDRE, R. (1998). Stochastic cohomology of the frame bundle of the loop space. *J. Nonlinear Math. Phys.* **5** 23–41. MR1609271

[40] LÉANDRE, R. (1998). Singular integral homology of the stochastic loop space. *Infin. Dimens. Anal. Quantum Probab. Relat. Top.* **1** 17–31. MR1611895

[41] LÉANDRE, R. (1999). Stochastic cohomology and Hochschild cohomology. In *Development of Infinite-Dimensional Noncommutative Analysis* (A. Hora, ed.) 17–26. RIMS, Kokyuroku. MR1752348

[42] LÉANDRE, R. (2000). Cover of the Brownian bridge and stochastic symplectic action. *Rev. Math. Phys.* **12** 91–137. MR1750777

[43] LÉANDRE, R. (2000). Anticipative Chen–Souriau cohomology and Hochschild cohomology. *Math. Phys. Stud.* **22** 185–199. MR1805914

[44] LÉANDRE, R. (2000). A sheaf theoretical approach to stochastic cohomology. *Rep. Math. Phys.* **46** 157–164. MR1803336

[45] LÉANDRE, R. (2001). Stochastic Adams theorem for a general compact manifold. *Rev. Math. Phys.* **13** 1095–1133. MR1853827

[46] LÉANDRE, R. (2001). Stochastic cohomology of Chen–Souriau and line bundle over the Brownian bridge. *Probab. Theory Related Fields* **120** 168–182. MR1841326

[47] LÉANDRE, R. (2001). Stochastic diffeology and homotopy. In *Stochastic Analysis and Mathematical Physics* (A. B. Cruzeiro and J.-C. Zambrini, eds.) 51–57. Birkhäuser, Boston. MR1886562

[48] LÉANDRE, R. (2002). Analysis over loop space and topology. *Math. Notes* **72** 212–229. MR1942549

[49] LÉANDRE, R. (2003). Stochastic algebraic de Rham complexes. *Acta Appl. Math.* **79** 217–247. MR2019477

[50] LÉANDRE, R. (2004). Hypoelliptic diffusion and cyclic cohomology. In *Stochastic Analysis* (R. Dalang, M. Dozzi and F. Russo, eds.) 165–185. Birkhäuser, Basel. MR2096288

[51] LÉANDRE, R. (2005). Speed of the Brownian loop on a manifold. *Quantum Probab.* (M. Bozejko, ed.). To appear.

[52] LÉANDRE, R. and SMOLYANOV, O. (1999). Stochastic homology of the loop space. In *Analysis on Infinite Dimensional Lie Groups and Algebras* (H. Heyer and J. Marion, eds.) 229–235. World Scientific, Singapore. MR1743171

[53] LODAY, J. L. (1998). *Cyclic Homology*, 2nd ed. Springer, New York. MR1600246

[54] NUALART, D. (1995). *The Malliavin Calculus and Related Topics.* Springer, New York. MR1344217

[55] NUALART, D. and PARDOUX, E. (1988). Stochastic calculus with anticipating integrands. *Probab. Theory Related Fields* **78** 535–581. MR950346

[56] RAMER, R. (1974). On the de Rham complex of finite codimensional forms on infinite dimensional manifolds. Thesis, Warwick Univ.

[57] SHIGEKAWA, I. (1986). De Rham–Hodge–Kodaira's decomposition on an abstract Wiener space. *J. Math. Kyoto Univ.* **26** 191–202. MR849215

[58] SMOLYANOV, O. G. (1986). De Rham current's and Stoke's formula in a Hilbert space. *Soviet. Math. Dokl.* **33** 140–144. MR834693

[59] SOURIAU, J. M. (1985). Un algorithme générateur de structures quantiques. *Astérisque* 341–399. MR837208

[60] SZABO, R. J. (2000). *Equivariant Cohomology and Localization of Paths Integrals in Physics. Lecture Notes in Phys.* **63**. Springer, Berlin. MR1762411

[61] WARNER, J. M. (1983). *Foundations of Differentiable Manifolds and Lie Groups.* Springer, New York. MR722297

[62] WITTEN, E. (1982). Supersymmetry and Morse theory. *J. Differential Geom.* **17** 661–692. MR683171

FACULTÉ DES SCIENCES
INSTITUT DE MATHÉMATIQUES DE BOURGOGNE
21000 DIJON
FRANCE
E-MAIL: Remi.leandre@u-bourgogne.fr